 \newcount\mgnf\newcount\tipi\newcount\tipoformule\newcount\greco
\tipi=2          
\tipoformule=0   
\global\newcount\numsec
\global\newcount\numfor
\global\newcount\numtheo
\global\advance\numtheo by 1

\def\senondefinito#1{\expandafter\ifx\csname#1\endcsname\relax}
\def\SIA #1,#2,#3 {\senondefinito{#1#2}%
\expandafter\xdef\csname #1#2\endcsname{#3}\else
\write16{???? ma #1,#2 e' gia' stato definito !!!!} \fi}
\def\etichetta(#1){(\veroparagrafo.\veraformula)%
\SIA e,#1,(\veroparagrafo.\veraformula) %
\global\advance\numfor by 1%
\write15{\string\FU (#1){\equ(#1)}}%
\write16{ EQ #1 ==> \equ(#1) }}

\def\letichetta(#1){\veroparagrafo.\verotheo
\SIA e,#1,{\veroparagrafo.\verotheo}
\global\advance\numtheo by 1
\write15{\string\FU (#1){\equ(#1)}}
\write16{ Sta \equ(#1) == #1 }}

\def\tetichetta(#1){\veroparagrafo.\veraformula 
\SIA e,#1,{(\veroparagrafo.\veraformula)}
\global\advance\numfor by 1
\write15{\string\FU (#1){\equ(#1)}}
\write16{ tag #1 ==> \equ(#1)}}

\def\FU(#1)#2{\SIA fu,#1,#2 }\def\etichettaa(#1){(A\veroparagrafo.\veraformula)%
\SIA e,#1,(A\veroparagrafo.\veraformula) %
\global\advance\numfor by 1%
\write15{\string\FU (#1){\equ(#1)}}%
\write16{ EQ #1 ==> \equ(#1) }}

\def\BOZZA{
\def\alato(##1){%
 {\rlap{\kern-\hsize\kern-1.4truecm{$\scriptstyle##1$}}}}%
\def\aolado(##1){%
 {
{
 \rlap{\kern-1.4truecm{$\scriptstyle##1$}}}}}
 }
\def\alato(#1){}
\def\aolado(#1){}
\def\veroparagrafo{\number\numsec}
\def\veraformula{\number\numfor}
\def\verotheo{\number\numtheo}

\def\Eq(#1){\eqno{\etichetta(#1)\alato(#1)}}
\def\eq(#1){\etichetta(#1)\alato(#1)}
\def\teq(#1){\tag{\aolado(#1)\tetichetta(#1)\alato(#1)}}
\def\Eqa(#1){\eqno{\etichettaa(#1)\alato(#1)}}
\def\eqa(#1){\etichettaa(#1)\alato(#1)}
\def\eqv(#1){\senondefinito{fu#1}$\clubsuit$#1
\write16{#1 non e' (ancora) definito}%
\else\csname fu#1\endcsname\fi}
\def\equ(#1){\senondefinito{e#1}\eqv(#1)\else\csname e#1\endcsname\fi}

\def\Lemma(#1){\aolado(#1)Lemma \letichetta(#1)}%
\def\Theorem(#1){{\aolado(#1)Theorem \letichetta(#1)}}%
\def\Proposition(#1){\aolado(#1){Proposition \letichetta(#1)}}%
\def\Corollary(#1){{\aolado(#1)Corollary \letichetta(#1)}}%
\def\Remark(#1){{\noindent\aolado(#1){\bf Remark \letichetta(#1).}}}%
\def\Definition(#1){{\noindent\aolado(#1){\bf Definition 
\letichetta(#1)$\!\!$\hskip-1.6truemm}}}
\def\Example(#1){\aolado(#1) Example \letichetta(#1)$\!\!$\hskip-1.6truemm}

\def\include#1{
\openin13=#1.aux \ifeof13 \relax \else
\input #1.aux \closein13 \fi}

\openin14=\jobname.aux \ifeof14 \relax \else
\input \jobname.aux \closein14 \fi
\openout15=\jobname.aux

{\count255=\time\divide\count255 by 60 \xdef\hourmin{\number\count255}
        \multiply\count255 by-60\advance\count255 by\time
   \xdef\hourmin{\hourmin:\ifnum\count255<10 0\fi\the\count255}}
\def\oramin{\hourmin }

\def\data{\number\day/\ifcase\month\or january \or february \or march \or april
\or may \or june \or july \or august \or september
\or october \or november \or december \fi/\number\year;\ \oramin}

\newcount\pgn \pgn=1
\def\foglio{\number\numsec:\number\pgn
\global\advance\pgn by 1}
\def\foglioa{A\number\numsec:\number\pgn
\global\advance\pgn by 1}

\footline={\rlap{\hbox{\copy200}}\hss\tenrm\folio\hss}

\def\TIPIO{
\font\setterm=amr7 
\def \settepunti{\def\rm{\fam0\setterm}
\textfont0=\setterm   
\normalbaselineskip=9pt\normalbaselines\rm }\let\nota=\settepunti}


\def\TIPITOT{
\font\twelverm=cmr12
\font\twelvei=cmmi12
\font\twelvesy=cmsy10 scaled\magstep1
\font\twelveex=cmex10 scaled\magstep1
\font\twelveit=cmti12
\font\twelvett=cmtt12
\font\twelvebf=cmbx12
\font\twelvesl=cmsl12
\font\ninerm=cmr9
\font\ninesy=cmsy9
\font\eightrm=cmr8
\font\eighti=cmmi8
\font\eightsy=cmsy8
\font\eightbf=cmbx8
\font\eighttt=cmtt8
\font\eightsl=cmsl8
\font\eightit=cmti8
\font\sixrm=cmr6
\font\sixbf=cmbx6
\font\sixi=cmmi6
\font\sixsy=cmsy6
\font\twelvetruecmr=cmr10 scaled\magstep1
\font\twelvetruecmsy=cmsy10 scaled\magstep1
\font\tentruecmr=cmr10
\font\tentruecmsy=cmsy10
\font\eighttruecmr=cmr8
\font\eighttruecmsy=cmsy8
\font\seventruecmr=cmr7
\font\seventruecmsy=cmsy7
\font\sixtruecmr=cmr6
\font\sixtruecmsy=cmsy6
\font\fivetruecmr=cmr5
\font\fivetruecmsy=cmsy5
\textfont\truecmr=\tentruecmr
\scriptfont\truecmr=\seventruecmr
\scriptscriptfont\truecmr=\fivetruecmr
\textfont\truecmsy=\tentruecmsy
\scriptfont\truecmsy=\seventruecmsy
\scriptscriptfont\truecmr=\fivetruecmr
\scriptscriptfont\truecmsy=\fivetruecmsy
\textfont\truecmr=\tentruecmr
\scriptfont\truecmr=\seventruecmr
\scriptscriptfont\truecmr=\fivetruecmr
\textfont\truecmsy=\tentruecmsy
\scriptfont\truecmsy=\seventruecmsy
\scriptscriptfont\truecmr=\fivetruecmr
\scriptscriptfont\truecmsy=\fivetruecmsy
\def \eightpoint{\def\rm{\fam0\eightrm}
\textfont0=\eightrm \scriptfont0=\sixrm \scriptscriptfont0=\fiverm
\textfont1=\eighti \scriptfont1=\sixi   \scriptscriptfont1=\fivei
\textfont2=\eightsy \scriptfont2=\sixsy   \scriptscriptfont2=\fivesy
\textfont3=\tenex \scriptfont3=\tenex   \scriptscriptfont3=\tenex
\textfont\itfam=\eightit  \def\it{\fam\itfam\eightit}%
\textfont\slfam=\eightsl  \def\sl{\fam\slfam\eightsl}%
\textfont\ttfam=\eighttt  \def\tt{\fam\ttfam\eighttt}%
\textfont\bffam=\eightbf  \scriptfont\bffam=\sixbf
\scriptscriptfont\bffam=\fivebf  \def\bf{\fam\bffam\eightbf}%
\tt \ttglue=.5em plus.25em minus.15em
\setbox\strutbox=\hbox{\vrule height7pt depth2pt width0pt}%
\normalbaselineskip=9pt
\let\sc=\sixrm  \let\big=\eightbig  \normalbaselines\rm
\textfont\truecmr=\eighttruecmr
\scriptfont\truecmr=\sixtruecmr
\scriptscriptfont\truecmr=\fivetruecmr
\textfont\truecmsy=\eighttruecmsy
\scriptfont\truecmsy=\sixtruecmsy }\let\nota=\eightpoint}

\newfam\msbfam   
\newfam\truecmr  
\newfam\truecmsy 
\newskip\ttglue
\ifnum\tipi=0\TIPIO \else\ifnum\tipi=1 \TIPI\else \TIPITOT\fi\fi


\def\sqr#1#2{{\vcenter{\vbox{\hrule height.#2pt
     \hbox{\vrule width.#2pt height#1pt \kern#1pt
   \vrule width.#2pt}\hrule height.#2pt}}}}
\def\qed{\hfill $\mathchoice\sqr64\sqr64\sqr{2.1}3\sqr{1.5}3$}



\newcount\foot
\foot=1
\def\note#1{\footnote{${}^{\number\foot}$}{\ftn #1}\advance\foot by 1}
\def\tag #1{\eqno{\hbox{\rm(#1)}}}
\def\frac#1#2{{#1\over #2}}

\def\text#1{\quad{\hbox{#1}}\quad}
\def\newpage{\vfill\eject}

\def\proof{{\noindent\pr Proof: }}

\def\thanks{\noindent{\bf Aknowledgements: }}
\font\pr=cmbxsl10


\font\ch=cmbx12
\font\ftn=cmr8

\font\it=cmti10
\font\bf=cmbx10
\font\sm=cmr7

%
\catcode`\X=12\catcode`\@=11
\def\n@wcount{\alloc@0\count\countdef\insc@unt}
\def\n@wwrite{\alloc@7\write\chardef\sixt@@n}
\def\n@wread{\alloc@6\read\chardef\sixt@@n}
\def\crossrefs#1{\ifx\alltgs#1\let\tr@ce=\alltgs\else\def\tr@ce{#1,}\fi
   \n@wwrite\cit@tionsout\openout\cit@tionsout=\jobname.cit 
   \write\cit@tionsout{\tr@ce}\expandafter\setfl@gs\tr@ce,}
\def\setfl@gs#1,{\def\@{#1}\ifx\@\empty\let\next=\relax
   \else\let\next=\setfl@gs\expandafter\xdef
   \csname#1tr@cetrue\endcsname{}\fi\next}
\newcount\sectno\sectno=0\newcount\subsectno\subsectno=0\def\r@s@t{\relax}
\def\resetall{\global\advance\sectno by 1\subsectno=0
  \gdef\firstpart{\number\sectno}\r@s@t}
\def\resetsub{\global\advance\subsectno by 1
   \gdef\firstpart{\number\sectno.\number\subsectno}\r@s@t}
\def\v@idline{\par}\def\firstpart{\number\sectno}
\def\l@c@l#1X{\firstpart.#1}\def\gl@b@l#1X{#1}\def\t@d@l#1X{{}}
\def\m@ketag#1#2{\expandafter\n@wcount\csname#2tagno\endcsname
     \csname#2tagno\endcsname=0\let\tail=\alltgs\xdef\alltgs{\tail#2,}%
  \ifx#1\l@c@l\let\tail=\r@s@t\xdef\r@s@t{\csname#2tagno\endcsname=0\tail}\fi
   \expandafter\gdef\csname#2cite\endcsname##1{\expandafter
     \ifx\csname#2tag##1\endcsname\relax?\else{\rm\csname#2tag##1\endcsname}\fi
    \expandafter\ifx\csname#2tr@cetrue\endcsname\relax\else
     \write\cit@tionsout{#2tag ##1 cited on page \folio.}\fi}%
   \expandafter\gdef\csname#2page\endcsname##1{\expandafter
     \ifx\csname#2page##1\endcsname\relax?\else\csname#2page##1\endcsname\fi
     \expandafter\ifx\csname#2tr@cetrue\endcsname\relax\else
     \write\cit@tionsout{#2tag ##1 cited on page \folio.}\fi}%
   \expandafter\gdef\csname#2tag\endcsname##1{\global\advance
     \csname#2tagno\endcsname by 1%
   \expandafter\ifx\csname#2check##1\endcsname\relax\else%
\fi
   \expandafter\xdef\csname#2check##1\endcsname{}%
   \expandafter\xdef\csname#2tag##1\endcsname
     {#1\number\csname#2tagno\endcsnameX}%
   \write\t@gsout{#2tag ##1 assigned number \csname#2tag##1\endcsname\space
      on page \number\count0.}%
   \csname#2tag##1\endcsname}}%
\def\m@kecs #1tag #2 assigned number #3 on page #4.%
   {\expandafter\gdef\csname#1tag#2\endcsname{#3}
   \expandafter\gdef\csname#1page#2\endcsname{#4}}
\def\re@der{\ifeof\t@gsin\let\next=\relax\else
    \read\t@gsin to\t@gline\ifx\t@gline\v@idline\else
    \expandafter\m@kecs \t@gline\fi\let \next=\re@der\fi\next}
\def\t@gs#1{\def\alltgs{}\m@ketag#1e\m@ketag#1s\m@ketag\t@d@l p
    \m@ketag\gl@b@l r \n@wread\t@gsin\openin\t@gsin=\jobname.tgs \re@der
    \closein\t@gsin\n@wwrite\t@gsout\openout\t@gsout=\jobname.tgs }
\outer\def\localtags{\t@gs\l@c@l}
\outer\def\globaltags{\t@gs\gl@b@l}
\outer\def\newlocaltag#1{\m@ketag\l@c@l{#1}}
\outer\def\newglobaltag#1{\m@ketag\gl@b@l{#1}}

\def\t@gsoff#1,{\def\@{#1}\ifx\@\empty\let\next=\relax\else\let\next=\t@gsoff
   \expandafter\gdef\csname#1cite\endcsname{\relax}
   \expandafter\gdef\csname#1page\endcsname##1{?}
   \expandafter\gdef\csname#1tag\endcsname{\relax}\fi\next}
\def\verbatimtags{\let\ift@gs=\iffalse\ifx\alltgs\relax\else
   \expandafter\t@gsoff\alltgs,\fi}
\catcode`\X=11 \catcode`\@=\active
\localtags
%
\setbox200\hbox{$\scriptscriptstyle \data $}
\global\newcount\numpunt
\magnification=1000
\hoffset=0.cm
\baselineskip=14pt  
\parindent=12pt
\lineskip=4pt\lineskiplimit=0.1pt
\parskip=0.1pt plus1pt

\hyphenation{small}
\newcount\mgnf\newcount\tipi\newcount\tipoformule\newcount\greco
\tipi=2          
\tipoformule=0   

\global\newcount\numsec
\global\newcount\numfor
\global\newcount\numtheo
\global\advance\numtheo by 1

\def\senondefinito#1{\expandafter\ifx\csname#1\endcsname\relax}
\def\SIA #1,#2,#3 {\senondefinito{#1#2}%
\expandafter\xdef\csname #1#2\endcsname{#3}\else
\write16{???? ma #1,#2 e' gia' stato definito !!!!} \fi}
\def\etichetta(#1){(\veroparagrafo.\veraformula)%
\SIA e,#1,(\veroparagrafo.\veraformula) %
\global\advance\numfor by 1%
\write15{\string\FU (#1){\equ(#1)}}%
\write16{ EQ #1 ==> \equ(#1) }}
\def\letichetta(#1){\veroparagrafo.\verotheo
\SIA e,#1,{\veroparagrafo.\verotheo}
\global\advance\numtheo by 1
\write15{\string\FU (#1){\equ(#1)}}
\write16{ Sta \equ(#1) == #1 }}
\def\tetichetta(#1){\veroparagrafo.\veraformula 
\SIA e,#1,{(\veroparagrafo.\veraformula)}
\global\advance\numfor by 1
\write15{\string\FU (#1){\equ(#1)}}
\write16{ tag #1 ==> \equ(#1)}}
\def\FU(#1)#2{\SIA fu,#1,#2 }
\def\etichettaa(#1){(A\veroparagrafo.\veraformula)%
\SIA e,#1,(A\veroparagrafo.\veraformula) %
\global\advance\numfor by 1%
\write15{\string\FU (#1){\equ(#1)}}%
\write16{ EQ #1 ==> \equ(#1) }}

\def\BOZZA{
\def\alato(##1){%
 {\rlap{\kern-\hsize\kern-1.4truecm{$\scriptstyle##1$}}}}%
\def\aolado(##1){%
 {
{
 \rlap{\kern-1.4truecm{$\scriptstyle##1$}}}}} 
}

\def\alato(#1){}
\def\aolado(#1){}

\def\veroparagrafo{\number\numsec}
\def\veraformula{\number\numfor}
\def\verotheo{\number\numtheo}

\def\Eq(#1){\eqno{\etichetta(#1)\alato(#1)}}
\def\eq(#1){\etichetta(#1)\alato(#1)}
\def\leq(#1){\leqno{\aolado(#1)\etichetta(#1)}}
\def\teq(#1){\tag{\aolado(#1)\tetichetta(#1)\alato(#1)}}
\def\Eqa(#1){\eqno{\etichettaa(#1)\alato(#1)}}
\def\eqa(#1){\etichettaa(#1)\alato(#1)}
\def\eqv(#1){\senondefinito{fu#1}$\clubsuit$#1
\write16{#1 non e' (ancora) definito}%
\else\csname fu#1\endcsname\fi}
\def\equ(#1){\senondefinito{e#1}\eqv(#1)\else\csname e#1\endcsname\fi}

\def\Lemma(#1){\aolado(#1)Lemma \letichetta(#1)}%
\def\Theorem(#1){{\aolado(#1)Theorem \letichetta(#1)}}%
\def\Proposition(#1){\aolado(#1){Proposition \letichetta(#1)}}%
\def\Corollary(#1){{\aolado(#1)Corollary \letichetta(#1)}}%
\def\Remark(#1){{\noindent\aolado(#1){\bf Remark \letichetta(#1).}}}%
\def\Definition(#1){{\noindent\aolado(#1){\bf Definition 
\letichetta(#1)$\!\!$\hskip-1.6truemm}}}
\def\Example(#1){\aolado(#1) Example \letichetta(#1)$\!\!$\hskip-1.6truemm}

\def\include#1{
\openin13=#1.aux \ifeof13 \relax \else
\input #1.aux \closein13 \fi}

\openin14=\jobname.aux \ifeof14 \relax \else
\input \jobname.aux \closein14 \fi
\openout15=\jobname.aux


{\count255=\time\divide\count255 by 60 \xdef\hourmin{\number\count255}
        \multiply\count255 by-60\advance\count255 by\time
   \xdef\hourmin{\hourmin:\ifnum\count255<10 0\fi\the\count255}}

\def\oramin{\hourmin }

\def\data{\number\day/\ifcase\month\or january \or february \or march \or april
\or may \or june \or july \or august \or september
\or october \or november \or december \fi/\number\year;\ \oramin}

\newcount\pgn \pgn=1
\def\foglio{\number\numsec:\number\pgn
\global\advance\pgn by 1}
\def\foglioa{A\number\numsec:\number\pgn
\global\advance\pgn by 1}

\footline={\rlap{\hbox{\copy200}}\hss\tenrm\folio\hss}

\def\TIPIO{
\font\setterm=amr7 
\def \settepunti{\def\rm{\fam0\setterm}
\textfont0=\setterm   
\normalbaselineskip=9pt\normalbaselines\rm }\let\nota=\settepunti}

\def\TIPITOT{
\font\twelverm=cmr12
\font\twelvei=cmmi12
\font\twelvesy=cmsy10 scaled\magstep1
\font\twelveex=cmex10 scaled\magstep1
\font\twelveit=cmti12
\font\twelvett=cmtt12
\font\twelvebf=cmbx12
\font\twelvesl=cmsl12
\font\ninerm=cmr9
\font\ninesy=cmsy9
\font\eightrm=cmr8
\font\eighti=cmmi8
\font\eightsy=cmsy8
\font\eightbf=cmbx8
\font\eighttt=cmtt8
\font\eightsl=cmsl8
\font\eightit=cmti8
\font\sixrm=cmr6
\font\sixbf=cmbx6
\font\sixi=cmmi6
\font\sixsy=cmsy6
\font\twelvetruecmr=cmr10 scaled\magstep1
\font\twelvetruecmsy=cmsy10 scaled\magstep1
\font\tentruecmr=cmr10
\font\tentruecmsy=cmsy10
\font\eighttruecmr=cmr8
\font\eighttruecmsy=cmsy8
\font\seventruecmr=cmr7
\font\seventruecmsy=cmsy7
\font\sixtruecmr=cmr6
\font\sixtruecmsy=cmsy6
\font\fivetruecmr=cmr5
\font\fivetruecmsy=cmsy5
\textfont\truecmr=\tentruecmr
\scriptfont\truecmr=\seventruecmr
\scriptscriptfont\truecmr=\fivetruecmr
\textfont\truecmsy=\tentruecmsy
\scriptfont\truecmsy=\seventruecmsy
\scriptscriptfont\truecmr=\fivetruecmr
\scriptscriptfont\truecmsy=\fivetruecmsy
\def \eightpoint{\def\rm{\fam0\eightrm}
\textfont0=\eightrm \scriptfont0=\sixrm \scriptscriptfont0=\fiverm
\textfont1=\eighti \scriptfont1=\sixi   \scriptscriptfont1=\fivei
\textfont2=\eightsy \scriptfont2=\sixsy   \scriptscriptfont2=\fivesy
\textfont3=\tenex \scriptfont3=\tenex   \scriptscriptfont3=\tenex
\textfont\itfam=\eightit  \def\it{\fam\itfam\eightit}%
\textfont\slfam=\eightsl  \def\sl{\fam\slfam\eightsl}%
\textfont\ttfam=\eighttt  \def\tt{\fam\ttfam\eighttt}%
\textfont\bffam=\eightbf  \scriptfont\bffam=\sixbf
\scriptscriptfont\bffam=\fivebf  \def\bf{\fam\bffam\eightbf}%
\tt \ttglue=.5em plus.25em minus.15em
\setbox\strutbox=\hbox{\vrule height7pt depth2pt width0pt}%
\normalbaselineskip=9pt
\let\sc=\sixrm  \let\big=\eightbig  \normalbaselines\rm
\textfont\truecmr=\eighttruecmr
\scriptfont\truecmr=\sixtruecmr
\scriptscriptfont\truecmr=\fivetruecmr
\textfont\truecmsy=\eighttruecmsy
\scriptfont\truecmsy=\sixtruecmsy }\let\nota=\eightpoint}

\newfam\msbfam   
\newfam\truecmr  
\newfam\truecmsy 
\newskip\ttglue
\ifnum\tipi=0\TIPIO \else\ifnum\tipi=1 \TIPI\else \TIPITOT\fi\fi

\def\chap #1#2{\line{\ch #1\hfill}\numsec=#2\numfor=1}


\def\AA{{\cal A}}
\def\BB{{\cal B}}
\def\CC{{\cal C}}
\def\DD{{\cal D}}
\def\EE{{\cal E}}

\def\GG{{\cal G}}

\def\LL{{\cal L}}

\def\SS{{\cal S}}

\def\NN{{\cal N}}
\def\MM{{\cal M}}

\def\LL{{\cal L}}

\def\a{\alpha}
\def\b{\beta}
\def\d{\delta}

\def\ve{\varepsilon}

\def\g{\gamma}

\def\l{\lambda}
\def\r{\rho}
\def\s{\sigma}
\def\t{\tau}
\def\th{\theta}

\def\D{\Delta}
\def\L{\Lambda}
\def\G{\Gamma}
\def\O{\Omega}
\def\S{\Sigma}

\def\h{\eta}

\def\E{{I\kern-.25em{E}}}
\def\N{{I\kern-.25em{N}}}
\def\M{{I\kern-.25em{M}}}
\def\R{{I\kern-.25em{R}}}
\def\Z{{Z\kern-.425em{Z}}}
\def\1{{1\kern-.25em\hbox{\rm I}}}
\def\eu{{1\kern-.25em\hbox{\sm I}}}
\def\C{{I\kern-.64em{C}}}
\def\P{{I\kern-.25em{P}}}
\def\T{{\rm I\kern -4.3pt{\rm T}}}
\def\W{{\rm W\kern -9.5pt{\rm W}}}
\def\J{{\rm I\kern -4.3pt{\rm J}}}
\def\B{{I\kern-.25em{B}}}

\def\X{{\rm X\kern -6.0pt{\rm X}}}
\def\V{{\rm V\kern -6.5pt{\rm V}}}
\def\U{{\rm U\kern -6.0pt{\rm U}}}
\def\Y{{\rm Y\kern -6.0pt{\rm Y}}}
\def\IG{\ {\rm I\kern -5.8pt{\rm G}}}
\def\IIG{{I\kern-.5em{G}}}


\def\rar{{\rightarrow}}

\def\leq{\le}
\def\geq{\ge}

\def\Es{\hbox{\bf E}}
\def\Pr{\hbox {\bf P}}


\def\rar{{\rightarrow}}

\def\leq{\le}
\def\geq{\ge}

\def\iL{\buildrel\,{\atop\circ}\over\L}
\def\Tr{\hbox{Tr}}

\def\mss{{  \nu_s}^{\a,N}}
\def\Qss{\text{\bf Q}\!\!\!\!\!\!}
\def\bro{\overline \rho}

\def\mss{{ \nu_s}^{\a,N}}
\def\Qss{\text{\bf Q}\!\!\!\!\!\!}
\def\bro{\overline \rho}
\hfuzz 16pt
\baselineskip 16pt

\catcode`\@=11  

\centerline {\bf   Lattice gas model in random medium and open boundaries:  } \centerline {\bf  hydrodynamic  and  relaxation to the steady state.    \footnote{$^*$}  
{\eightrm work supported by   
   INDAM-CNRS, Roma TRE,  University of Rouen. }}
\vskip1cm 
\centerline{ 
Mustapha  Mourragui     \footnote{$^1$}{\eightrm 
Universit\'e de Rouen, LMRS, UMR 6085,  Avenue de l'Universit\'e, 
BP. 12,  76801, Saint Etienne du Rouvray, France.}
\footnote{}{\eightrm Mustapha.Mourragui@univ-rouen.fr}
\hskip.2cm 
 and Enza Orlandi  \footnote{$^2$}{\eightrm 
Dipartimento di Matematica, Universit\'a di Roma Tre,  L.go S.Murialdo 1, 
00146 Roma, Italy.    orlandi@mat.uniroma3.it }
 }
\footnote{}{\eightrm {\eightit Key Words }:
Random environment,  Nongradient systems,   Stationary nonequilibrium states. 
 }
\footnote{}{\eightrm {\eightit 2000 Mathematics Subject Classification}. Primary 82C22,
60K35,    Secondary 60F10, 82C35.}
\vskip.5cm

{ \bf Abstract} 
  We consider a lattice gas interacting by the exclusion
rule in the presence of a  random field given by 
i.i.d. bounded random variables  in a bounded domain
in contact with particles    reservoir   at  different densities.   We show, in dimensions  $d \ge 3$,    that   the rescaled  empirical density
field  almost surely,  with respect to the random field,      converges    to the unique weak solution of a non linear  parabolic    equation    having  the diffusion matrix 
determined  by the statistical properties of the external random field and
   boundary conditions determined by the density of the reservoir.      Further we 
show that 
the rescaled  empirical density field,  in the stationary regime, 
 almost surely   with respect to the random field, converges to the solution of  the  associated stationary transport equation.   
\bigskip \bigskip  
\chap {1 Introduction}1
\numsec= 1 \numfor= 1 \numtheo=1
In the last years there has been  several papers devoted  in  understanding macroscopic properties of 
non equilibrium systems. Typical examples are  systems   in contact   with two thermostats at  different temperature or   with two reservoirs  at different  densities.   
   A mathematical model of open systems is provided by stochastic models of interacting particles systems  performing  a  local reversible dynamics  (for example a reversible hopping dynamics) in a domain and some external mechanism of creation and annihilation of particles on the boundary of the domain, modeling the reservoirs,  which makes the full process non reversible.     The first question that one might ask for these systems  is  the  derivation of  the hydrodynamic behavior   (law of large number)  for the locally conserved field in the non stationary and stationary regime.       There has been  important classes of models, see for example [ELS1,2] , [DFIP], [KLO]  in which it has been proved  the law of large numbers for the empirical density in the   stationary  regime.  Typical generic  feature
    of these systems  is that they exhibit  long range correlation in their steady state.    These long range correlations have been calculated from the microscopic dynamics only in very few cases;
   mainly in the case of the symmetric exclusion process,  see  [Sp1],  the  asymmetric exclusion process, see [DEL],
   and in the weakly   asymmetric exclusion process, see [DELO].   
   More recently  breakthroughs  were achieved analyzing the large deviations principle for the stationary measure. 
   We refer to [BSGJL] for a review of  works on the statistical mechanics of non equilibrium processes based on the analysis of large deviations properties of microscopic systems.
  
 In this paper we focus on the first step.   We derive the  macroscopic limit  in the stationary and not stationary regime (hydrodynamic   limit)
for  a  particles system   evolving according to  local- conservative dynamics (Kawasaki)  with hard core exclusion rule  and with rates depending  on a quenched random field  in a cylinder domain $d\ge 3$ in which the basis, denoted $\Gamma$, are kept at different densities.  The restriction on the  dimensions is  only   technical. We comment on this later. The rates are chosen so that the system satisfies a detailed balance condition with respect to a family of random Bernoulli measures (the random field Ising model at infinite temperature).  To model  the presence of the reservoirs, 
as in previous papers,  we  superimpose  at  the boundary,  to the local-conservative dynamics, a jump dynamics (creation and destruction of particle). The rates of the birth and death process    depend on   the  realizations of the random field   and are chosen so that   a   random Bernoulli measure with a suitable choice of the  chemical potential  is  reversible for  it.  This latter  dynamic is of course not conservative and  keeps the fixed value of the  density on  the boundary.  There is a flow of  density through the  full system and the full  dynamic  is not   reversible.  
The   bulk dynamic      models electron transport in doped crystals. In this case the exclusion rule is given by the Pauli principle and the presence of impurities in the crystals is the  origin of the presence of  quenched random field,  see [KW].
The presence of the  random field together with the exclusion rule makes the problem high not trivial. The transport properties of such systems in the case  of periodic boundary condition on $\Gamma$ has been studied by Faggionato  and Martinelli,  [FM].   They derived in $d \ge 3$,  the   hydrodynamic   limit  and  gave a variational formula  for  the bulk diffusion,  equivalent to the Green-Kubo formula. They  proved  that the bulk diffusion is a deterministic quantity depending    on the statistical properties of the   random field.    
  Later,  Quastel [Q] derived   in all dimensions for the same model investigated by   [FM]  the hydrodynamic   limit for the local empirical density.   Applying    the  method proposed by Quastel,  we  could  extend      our  results  in all dimensions.  Since our aim is to understand the role of  the randomness 
 in the non stationary  and stationary  state  and not the role of dimensions in the bulk dynamics  we  state and prove our results  in $d \ge 3$. 
 Dynamical Large deviations  for the  same model
and always with periodic boundary conditions have been derived in [MO] as special case of a   more general system  discussed there. 
 The  bulk dynamics  is    of the so-called  {\sl nongradient} type.   
Roughly speaking, the gradient condition says that the  microscopic current 
is already the gradient of a function of the density field. 
Further  it is not   translation invariant, for a given
disorder 
configuration.      In order to prove the hydrodynamic behavior of the system, we follow the entropy method introduced by Guo, Papanicolaou and Varadhan [GPV].  It relies on an estimate of the entropy of the states of process  with respect to a reference invariant state.  By the general theory of Markov Processes the entropy of the state of a process with respect to an invariant state decreases in time.
The main problem is that in the model considered the reference invariant state is not explicitly known.  To overcome this difficulty  we  compute the   entropy of the state of the process with respect to a 
 product measure with slowly varying profile.   Since this measure is not invariant, the entropy does not need to decrease  and   we need to  estimate the rate at which it increases.     This type of strategy has been used in previous papers dealing with the same type of problems,   see [KLO] and [LMS],  which considered generalized exclusion process of non gradient type.  
 The main difference with the previous mentioned papers is the presence of the randomness in the model    considered here.   This forces to  take   on the boundary  a jump process   depending on the external random field. 
 Important step to derive the final  results  is then   a   convenient  application of the    ergodic theorem.   Further we show that the empirical density field obeys a law of large numbers with respect to the  stationary random measures  (hydrostatic).   This is achieved  proving   that  it is possible to  derive the hydrodynamic for  the  evolution of the  empirical  measures starting from any initial particle configurations distributed according to the   stationary measure, even  though   it  is not possible to  identify the      profile.  Then we exploit that the stationary solution of the parabolic nonlinear equation is unique and is a  global attractor for the  evolution.   These two ingredients allow
 to conclude.  Similar strategy for  proving the  hydrostatic is   used  in the paper in preparation  by  
 Farfan Vargas,  Landim and  Mourragui,  [FLM].   
 
\bigskip

\bigskip
\bigskip


\chap{2 The model  and the main results}2
\numsec= 2
\numfor= 1
\numtheo=1

\bigskip

\noindent{\it 2.1. The model}
\medskip

We consider the $d-$ dimensional lattice $ \Z^d$  with sites $x= (x_1,\dots, x_d)$ and 
canonical basis   $\EE=\{e_1,\ldots ,e_d\}$ and we assume  in all the paper that   $d\ge 3$.  We
denote by $\L:= [-1,1 ]\times \T^{d-1}$, where $\T^{d-1}$ is
the  $(d-1)$-dimensional torus  of  diameter 1 and by  $\Gamma$ the
boundary of $\L$.

Fix an integer   $ N \ge 1 $.   Denote by $\L_N\equiv \{-N,\cdots,N\}\times \T_N^{d-1}$ the
cylinder in $\Z^{d}$ of length $2N+1$ with basis the $(d-1)$-dimensional
discrete torus $\T_N^{d-1}$   
and by $\G_N=\{ x\in \L_N\,|\, x_1=\pm N\}$ 
the boundary of $\L_N$.  The elements of $\L_N$ will
be denoted by letters $x,y,\ldots$ and the  elements of $\L$ by $u, v, \ldots$.

For a fixed $A> 0$, let $\Sigma_D = [-A,A]^{\Z^d}$ 
be  the set of disorder configurations on $\Z^d$.
On $\Sigma_D$ we define a product, translation invariant probability 
measure $\P$. We denote by $\E$ the expectation with respect to $\P$, 
and by $\a\equiv\{\a(x),\ x\in \Z^d\}$, $ \a(x) \in [-A,A]$, a  disorder configuration
in $\Sigma_D$.
A configuration $\a\in\Sigma_D$ induces in a natural way a disorder 
configuration   $ \a_N$ on $ \L_N$, by identifying a cube  centered 
at the origin of side $2N+1$ with $\L_N$. 
By a slight abuse  of notation whenever in the following  we refer to a disorder
configuration   either on $\L_N$ or on $\Z^d$  we denote it  by $\a$.  
We denote by $\SS_N \equiv \{0, 1\}^{\L_N}$ and  
$\SS\equiv \{0, 1\}^{{\Z^d}}$  the configuration spaces, both   
equipped with the product topology; elements of   $\SS_N$ or  $\SS $ are denoted by $\eta$, so that $\eta(x)=1$, resp $0$, if the site $x$ is occupied, resp empty, for the configuration $\eta$.  
Given  $\a\in \Sigma_{D}$, we consider the random  Hamiltonian $H^\a: \SS_N \to \R$,  
$$
H^\a(\h)
=- \sum_{x\in \L_N} \a(x)\h (x). 
\Eq(2.1a)
$$  
We denote by  
 $\mu^{\a,\l}_{N}$ the grand canonical random  Gibbs measure on
 $\SS_N$ 
associated to the Hamiltonian \eqv (2.1a) with chemical potential $\l\in \R$, i.e the  random  Bernoulli product 
measure 
  $$
\mu^{\a,\l}_{N}(\h)
=  \prod_{ x \in \L_N}   \Big\{
 \frac{ e^{ [\a (x) + \l ] \eta (x)}  } {  
e^{ [\a (x) + \l ] } +1   } \Big\}.  \Eq(2.4)
$$
When $\l=0$, we simply write $\mu^\a_N$.  We denote   by $\mu^{\a,\l}(\cdot) $ and when  $ \l=0$,   
  $\mu^\a (\cdot) $    the measure \eqv (2.4)
 on the infinite product space $\SS$.  Moreover, for a probability measure $\mu$ and a bounded 
function $f$, both defined on  $\SS$ or $\SS_N$, we denote by 
$\Es^{\mu}(f)$ the expectation  of $f$ with respect to $\mu$. 
We need  to introduce also the canonical measures $\nu_{\rho}^{\a,N}$,
$$ \nu_{\rho}^{\a,N}(\cdot) = \mu^{\a,\l}_{N}(\cdot | \sum_{x \in \L_N} \eta_x = \r |\L_N| ) $$ 
for $ \r \in [0, \frac 1 {|\L_N|},\dots,1]$.  It is well known [CM]  that the canonical and the grand canonical measures are closely related if the chemical potential $\l$  is  chosen canonical conjugate  to the density $\r$,  in the sense that the average density with respect to
$\mu^{\a,\l}_{N}$ is equal to $\r$.   So as in [FM] one can define the random empirical chemical potential and
the {\sl annealed chemical potential} $\l_0(\r)$. To our aim it is enough to consider $\l_0(\r)$. 
For 
$\r\in[0,1]$, the function $\l_0(\r)$ is defined     as the unique $\l$ so that   
$$ 
\E\left[\int\h(0)d\mu^{\a,\l }(\h)\right] 
 =\E\left[\frac{e^{\a(0)+\l }}{1+e^{\a(0)+\l }}\right] 
=\r\; .\Eq (pot1) 
$$
We will consider as reference measure  the random Bernoulli product
measure  $\nu_{\rho(\cdot)}^{\a,N}$ on $\SS_N$   defined  for positive profile 
$\rho:\L\to (0,1)$ by 
$$
\nu_{\rho(\cdot)}^{\a,N} (\eta) = \prod_{ x \in \L_N}   \Big\{
 \frac{ e^{ [\a (x) + \l_0(\rho(x/N)) ] \eta (x)}  } {  
e^{ [\a (x) + \l_0(\rho(x/N)) ] } +1   } \Big\}\; ,
\Eq (inv1)$$
if $\rho(\cdot) \equiv \rho$ is constant, we shall denote simply
$
\nu_{\rho(\cdot)}^{\a,N}=\nu_{\rho}^{\a,N}\; .
$
We denote by    $\h^{x,y}$  the configuration obtained from $\h$ by interchanging 
the values at $x$ and $y$:
$$
 \h^{x,y} (z) = \left \{ \eqalign 
{ & \h (x) \qquad \hbox { if } \qquad  z=y \cr 
&  \h (y) \qquad \hbox { if }\qquad  z=x \cr  
&  \h (z) \qquad \hbox { otherwise, }\; 
  }\right. \Eq (C.1a)
  $$
and by 
$\h^{x}$   the configuration obtained from
$\h$ by flipping the occupation number at site $x$:
$$
\h^x(z)= \left \{ \eqalign 
 { &\h (z)  \qquad \hbox {if}\qquad z\ne x \cr
 & 1-\h(x) \qquad \hbox {if}\qquad z=x .
    } \right. \Eq (C.1b)
$$
 Further, for $f:\SS_N \to\R$, $x,y\in\L_N$,   we denote
   $$
(\nabla_{x,y}f)(\h) = f(\h^{x,y})-f(\h). 
$$
The disordered  exclusion process on $\L_N$ 
with random  reservoirs at its boundary $\Gamma_N$ is 
the  Markov process on $\SS_N$ whose generator  $\LL_N$ 
can be decomposed  as 
$$
\LL_N=\LL^0_{N}+\LL^b_{N}\; ,\Eq (gener1)
$$
where the generators $\LL^0_{N}$,  $\LL^b_{N}$ act on function $f:\SS_N \to\R$
as 
$$
\left(\LL^0_{N}  f\right) (\h)= 
 \sum_{e\in\EE}\sum_{x \in \L_N, x+e \in \L_N}  C_N (
x,x+e; \h) \left[ (\nabla_{x,x+e}f)(\h) \right] \; ,\Eq (2.5)   
$$
where $e$ is  a generic element of $\EE$, the 
   rate $C_N $ is given by 
$$   
C_N (x,y;\h)\equiv C_N^\alpha (x,y;\h) = 
\exp\Big\{- \frac{(\nabla_{x,y}H^\a )(\h)}{2} \Big\}\; ; \Eq (rate) $$
and
$$
\left(\LL^b_{N} f\right)(\eta) =
    \sum_{x\in\G_N} C^{b}( x/N,\eta) \big[f(\eta^x)
    -f(\eta)\big]\, . \Eq (2.51)
$$
To define the rate $C^{b}( x/N,\eta)$ we 
 fix a   function $b(\cdot)$   on $\Gamma$, representing the density of  the reservoirs.  We   assume that $b (\cdot)$  is  the restriction on $\G$ of a  smooth function  $\gamma (\cdot)$ defined  on a neighborhood $V$ of $\L$,    $ \g : V \to  (0,1)$   and  
  $\gamma (u) = b(u) $  for   $u \in \Gamma$.
The rate $C^{b}$ is  chosen so that  $\LL^b_{N}$ is reversible with
respect to     $\nu_{\gamma(\cdot)}^{\a,N}$   
$$
C^{b} (x/N, \h) =
\h (x) \exp\Big\{ -\frac{\a (x) +\l_0(b (\frac x N)) }{2}\Big\} +(1-\h (x)) 
 \exp\Big\{ \frac{\a(x) +\l_0(b(\frac x N) )}{2}\Big\} \, . \Eq (CD1)
$$
The first term in \eqv (CD1) is the  creation   rate, the second one is the annihilation  rate. 
Next we recall the  relevant properties of 
  $C_N (x,y;\h)$:

a) {\it detailed balance condition with respect to the measure  \eqv (2.4)},

b) {\it positivity and boundedness}: there exists $a>0$ such that 
 $$
a^{-1} \le C_N (x,y;\h) \leq a, \Eq (phi-pos)
  $$

c) {\it translation covariant:}  
$$ 
C_N^{\a}(x,y;\h) 
=  C_N^{\t_z\a}(x-z,y-z; \t_z\h) 
= \t_z C_N^{\a}(x-z,y-z;\h)\; , \Eq (C.3)
$$ 
where for  $z$ in $\Z^d$,  $\t_z$ denotes the space shift by $z$
units on $\SS\times \Sigma_D$ defined for all $\h\in\SS$, $\a\in \Sigma_D $ and
$g:\SS\times\Sigma_D\to\R$  by
$$
(\t_z\h)(x)=\h(x+z), \; (\t_z\a)(x)=\a(x+z), \; (\t_z g)(\h,\a)=
g(\tau_z \h,\tau_z\a) \; .\Eq (shift) 
  $$
  We omit  to write in the notation the explicit dependence on the randomness $\a$, unless
there is an ambiguity. 
   The  process  arising from the full generator \eqv (gener1) is then a superposition of a dynamics with a conservation law (the Kawasaki random dynamics) acting on the whole $\L_N$ and
a  birth and death  process  acting on   $\Gamma$.
Remark that if $b (\cdot) \equiv b_0$ for some positive  constant $b_0$, then
the    generator $\LL_{N} $, see \eqv (gener1),  is self-adjoint in 
$L^2(\nu_{b_0}^{\a,N} )$ and  the measure $\nu_{b_0}^{\a,N}$
is the  stationary measure  for the full dynamics $\LL_N$.  In the general case, when $b (\cdot)$ is not constant,   since the Markov process on  $\SS_N$  with  generator  \eqv (gener1),  is  irreducible for all $N\ge 1$, there exists always an unique invariant measure but in general cannot be written in an explicit form.

 \vskip0.5cm 

\medskip
 \medskip
\noindent{\it 2.2. The macroscopic equation }
\medskip
The macroscopic evolution of  the  local particles density $\r $  is described by 
 the  quasi linear parabolic equation   
 $$
 \left\{ \eqalign{
&\partial_t \rho \; = \; \nabla\cdot\Big( D(\rho) \nabla \rho \Big)\, , \cr
&\rho (0 ,\cdot) \; =\; \rho_0\, ,   \cr
& \rho (t, \cdot){\big\vert_\Gamma} \; =\; b (\cdot) \quad \hbox{for}
  \quad  t> 0 \;, }
     \right. 
\Eq(heq1)
$$
 where  $D(\rho) $ is  the  diffusion  matrix  given  in \eqv (dif),   $b (\cdot)\in C^1(\G)$  represents  the interaction  with the reservoirs appearing  as boundary conditions to be imposed on the solution, see  its    definition before \eqv (CD1),  and    
$\rho_0:   \Lambda \to  [0,1] $ is  the  initial profile. 
 The diffusion matrix is the one derived in  
 [FM].      To define it,  let \footnote{$^*$}   { \eightrm  A function $g: \SS\times\L_D\to\R $ is local if the support of $g$, $\D_g
$, i.e. the smallest subset of $ \Z^d$ such that $g$ depends only on $ \{ (  \eta
(x),\a(x) ) 
\; x \in \D_g \} $, is  finite. The function $g$ is bounded if $ \sup_\eta \sup_\a
|g(\eta,\a)| < \infty $. }  
  $$ 
\IG\equiv\{ g: \SS\times\L_D\to\R;\ \ \hbox {local and bounded }\} \;,\Eq(3.13)
$$
and  for $g
\in
\IG
$, $
\G_g(\h)=\sum_{x\in {\Z}^d} \big(\t_x g\big)(\h,\a)
$.   
The  $
\G_g(\h)$ is  a formal expression, but the  difference 
$\nabla_{0,e}\G_g(\h) =\G_g(\h^{0,e}) -\G_g(\h) $  for $e\in\EE$  
is  meaningful. For each $\r\in(0,1)$, let $D(\r) =\{ D_{i,j}(\r),\ 
1\le i,j\le d\}$ be the symmetric matrix defined, for every $a\in\R^d$,  
by the variational formula
$$ 
(a\cdot D(\r)a)= \frac 1 {2\chi(\r)} \inf_{g \in\IIG }  \sum_{i=1}^d
\E \left [ 
\Es^{\mu^{\a,\l_0(\r) }}\left(
C^{0} (0,e_i;\h)\Big\{ a_i \nabla_{0,e_i}\h(0)+ 
   (\nabla_{0,e_i} \G_g)(\h)\Big\}^2 \right ) \right ]\Eq(dif)
$$
where
$\l_0(\r)$ is defined in \eqv (pot1), 
  $\chi(\rho)$ is the static compressibility   given by
 $$
 \chi(\r)=\E\left[
\int \h(0)^2d\mu^{\a, \l_0(\r)}(\h)
-\left (\int \h(0)d\mu^{\a, \l_0(\r)}(\h)\right)^2
\right]\;,       \Eq (comp) 
$$
for $a,b\in\R^d$,
$(a \cdot b) $ is the  scalar vector product of $a$ and $b$
  and, recall, $\Es^{ \mu^{\a,\l_0 (\r) } }(\cdot)$ is the expectation with respect to  
$\mu^{\a,\l_0(\r)}$, see after \equ(2.4), 
the random Bernoulli product measure  on $\SS$ with annealed  chemical potential $\l_0(\r)$.  In  Theorem 2.1 of [FM] 
it has been proved,  for
$d\ge 3$ and for
$\r
\in (0,1)
$, the existence of the symmetric diffusion matrix defined in \eqv (dif).   Further it has been 
proved that   the 
coefficients
$D_{i,j}(\cdot)$ are  nonlinear continuous functions  in the open interval  $(0,1) $ and 
there exists a   constant $C>1$, depending on dimensions and bound on the 
random field, such that 
$$
\frac {  \1 } { C } \le D(\r) \le C  \1  \qquad \r \in (0,1) \Eq (elliptic) 
$$ 
where $\1$ is the $d\times d$ identity matrix.  One expects  the matrix  $D (\cdot)$ to be  extended continuously
 to the closed interval $[0,1]$ and actually  to   be  
  a smooth function of   $\r$, [KW].  We will assume  all trough the paper that $D (\cdot)$ is well defined  in $[0,1]$ and Lipschitz in the open interval. 
  The  diffusion matrix
$D(\r)$ in a   solid, in a regime of linear response,   is
linked  to the  mobility   
$\frac 12 \s(\r)$, see [Sp],  via the Einstein relation 
$$D(\r)= \frac 12 \s(\r) \chi(\r)^{-1}.  \Eq (sp) $$
 The   $\chi(\rho)$ is a  smooth function of $\r$ in $[0,1]$  and it
 can be easily proven  from
\eqv (comp) that   
$$ \frac 1 2   \r (1-\r) \le  \chi (\r) \le    \r (1-\r); \qquad  \frac 1 C \r (1-\r)  \1 \le \s(\r)\le C  \1 \r
(1-\r), 
\Eq (rouen7a)
$$  
 where      $C $  is   a constant that may change from one occurrence to the next. 

\vskip1.cm  {\it Weak solutions} By weak solution of  \eqv (heq1) we mean a 
 function $\rho(\cdot, \cdot ):[0,T]\times \L\to \R $ satisfying
 \smallskip
\noindent {(IB1)} $\rho \in L^2 \left(( 0,T) ;H^1(\L)\right)$ :
$$
\int_0^T d s\Big( \int_\L {\parallel\nabla \rho(s,u)\parallel}^2 
du \Big)<\infty \; ; \Eq(energy)$$

\smallskip
\noindent {(IB2)}   For every function $G(t,u)=G_t(u)$ in 
${\CC}_c^{1,2}\big([0,T]\times\iL \big)$, where $\iL=]-1,1[ \times \T^{d-1}$ and ${\CC}_c^{1,2}\big([0,T]\times \iL \big)$ is 
the space  of functions  from 
  $[0,T]\times \iL$ to $\R$  twice continuously differentiable in $\L$
with continuous time derivative and having compact support in $\iL$ we have

$$\eqalign{
& \int_\L du \big\{ G_T(u)\rho(T,u)-G_0(u)\rho(0,u)\big\} -
\int_0^T ds \int_\L d u \,  (\partial_s G_s)(u)\rho(s,u) \cr
& \quad =\; -\int_0^T d s \Big\{\int_\L d u \, D(\rho(s,u))\nabla\rho(s,u)\cdot
\nabla G_s(u)  \Big\}\ ;
}$$
  \smallskip
\noindent {(IB3)} For any $t\in  (0,T]$,
 $\Tr(\rho (t,\cdot))=b(\cdot),\ \hbox{\rm a.e.}$.
  \smallskip
\noindent {(IB4)}  $\rho (0,u)= \rho_0(u)$. a.e.
 \smallskip
 
 Notice that, since the original particle model cannot have more than one particle at a
lattice site any  solution $\r $ of \eqv (heq1) is bounded between 0 and 1.
  The existence and uniqueness of the weak solution of \eqv  (heq1) when    
\eqv (elliptic) holds and $D(\cdot)$ is 
Lipschitz continuous for $\r \in (0,1)$,     can be done using standard
analysis tools.  We refer to [LSU], chapter V  or [DL]. 
Further, one  immediately  obtains    by the characterization of $H^{-1} (\L)$, see for example [E], page 283,  that
 $\partial_t \rho \in L^2\big(0,T;H^{-1}(\Lambda)\big)$. 
Recall that $H^{-1}(\Lambda)$ is  the dual of $H^1_0(\Lambda)$, i.e.  the  Banach space
equipped with the norm
$$
\|v\|_{-1} \;=\; \sup_{f}
\Big\{  \big<v,f\big>:  \|f \|_{H^1_0(\L)} \le 1 \Big \}\;.
\Eq(norm-1) 
$$
  \medskip
\noindent  {\it  Stationary  solution}  We denote by $ \bar \r$ the stationary solution of  \eqv (heq1), i.e.  a function from  $\L \to [0,1]$  so that   
$\bar \rho \in H^1(\L)$,  for 
 $G \in   \CC^{2}_c (\iL \big )$ 
we have 
$$  \left \{ \eqalign { &   \int_\L d u \, D(\bar \rho(u))\nabla {\bar \rho}( u)\cdot
\nabla G (u)   \; =0,\cr &  
 \Tr(\bar \rho (\cdot))=b(\cdot),\ \hbox{\rm a.e.} } \right. \Eq (SS.1)$$

 \medskip
\noindent{\it 2.3. The main results }
\medskip
 For any $T>0$,  we denote by  $(\h_t)_{t\in [0,T]}$ 
the Markov process  on $\SS_N $ with generator 
$N^{2}\LL_N $ starting from $\eta_0= \eta$ and     by $\Pr_{\eta }:=\Pr^\a_{\eta } $  its distribution when the initial configuration is $\eta$.   We  remind that we omit to write explicitly the dependence on $\a$.   The $\Pr_{\eta } $ is a 
probability measure on the path space $D( [0,T],\SS_N)$, which we consider endowed with  the  Skorohod topology and the corresponding Borel $\s-$algebra. Expectation with respect to  $\Pr_{\eta } $ is denoted by $\Es_{\eta } $. 
   If $ \mu^N$ is a probability measure on $ \SS_N$ we denote   $\Pr_{\mu^N} (\cdot) = \int_{\SS_N} \Pr_{\eta }( \cdot) \mu^N (d\eta)   $
  and by 
$\Es_{\mu^N } $ the expectation with respect to 
$\Pr_{\mu^N } $.  For $t\in[0,T],\h\in \SS_N$, let the empirical 
measure $\pi_t^N$ be defined by 
 $$ 
\pi_t^N(\eta) \equiv \pi^N(du;\h_t)\,
=\, \frac{1}{N^d}\sum_{x\in{\L_N }}\h_t(x)\, \d_{ x/N} (d u)\;, \Eq (E.5)
$$ 
where  
$\d_{u}(\cdot)$ is the Dirac measure on $\L$ concentrated on $u$.
Since $\h(x)\in\{0,1\}$, relation  \eqv (E.5) induces from  
$\Pr_{\mu^N} $ a distribution $Q_{\mu^N} $ on the Skorohod space  
$D([0,T],\MM_1(\L))$, where  $\MM_1(\L)$ is the set  of positive Borel measures    on $\L$
with total mass  bounded by 1, endowed  with the weak topology.
Denote by  $ \MM^0_1(\L)$ the subset of $\MM_1(\L)$ of all absolutely continuous measures w.r.t. the Lebesgue
measure with density bounded by 1:
$$ \MM^0_1(\L)= \left \{ \pi \in \MM_1(\L): \pi(du)= \r(u) du \quad \hbox {and } \quad 0\le \r(u) \le 1
\quad \hbox {a.e. } \ \right \}\; , $$
$ \MM^0_1(\L) $ is a closed subset of  $\MM_1(\L)$  endowed with the weak topology and 
$D([0,T],\MM^0_1(\L))$ is a closed subset of $D([0,T],\MM_1(\L))$ for the Skorohod topology. 
 To state next theorem we need the following definition. 
\vskip0.5cm \noindent

{\bf Definition } { \it  Given a Lebesgue absolutely continuous measure $ \r (u) du \in \MM^0_1(\L) $, 
 a  sequence of probability measures $(\mu^N)_{N \ge 0}$ on $\SS_N$  is said to correspond to the
macroscopic profile $\r$ if,
  for  any smooth  function $G$ and  
$\d>0$ }
$$
\lim_{N \to \infty} \mu^N
\Big\{ \Big| \frac 1 {N^d} \sum_{ x \in \L_N}  G( x/N) \eta (x)  -\int_\L G(u) \r (u) du  \Big|>\d \Big\}=0.  \Eq (o.1)
$$

\medskip
\noindent{\bf  \Theorem (th-hydro)} {\it Let $d\ge 3$
and assume that $D(\r)$ can be continuously extended to the
closed interval $[0,1]$.  Let $\mu^N$ be  a sequence of probability measures  on
$\SS_N$  corresponding  to the initial profile $\rho_0$.
Then, $\P$ a.s. 
the sequence of probability measures $({Q_{\mu^N} } )_{N\geq 0}$ 
is tight and all its limit points $Q^*$ are concentrated 
on  $\r(t,u){\rm d}u$, whose densities
are weak solutions of the equation \equ(heq1).
Moreover if   $D(\cdot) $ is 
Lipschitz continuous for $\r \in  (0,1)$, then $(Q_{\mu^N} )_{N\geq 0}$
converges weakly,   as $N\uparrow \infty $, to $Q^*$. This limit point is concentrated on 
the unique weak solution of equation \equ(heq1).  
\medskip
    }

Denote by $ \nu_{s}^{\a, N} $  the unique invariant measure of the Markov process $(\h_t)_{t\in [0,T]}$ with generator 
$N^{2}\LL_N $. We have the following: 
\medskip
\noindent{\bf  \Theorem (hydro1)}  {\it     Let $d\ge 3$,
 assume that $D(\r)$ can be continuously extended to the
closed interval $[0,1]$ and    
Lipschitz continuous for $\r \in (0,1)$. For every continuous function $G: \L \to \R$ and every  $ \d>0$,
 $$
\lim_{N \to \infty} \nu_{s}^{\a, N}
\Big\{ \Big| \frac 1 {N^d}  \sum_{ x \in \L_N}  G( x/N) \eta (x)  -\int_\L G(u) \bar  \r (u) du  \Big|>\d \Big\}=0, \quad \P=1,   \Eq (s.1)
$$
with  $\bar  \r (\cdot) $  satisfying \eqv (SS.1).  } 
\bigskip
\bigskip


\chap{3. Strategy of proof and basic estimates }3
\numsec=3
\numfor=1
\numtheo=1

\bigskip

\noindent{\sl 3.1. The steps to prove Theorem \equ(th-hydro)}
\medskip 
 To prove the hydrodynamic behavior of the system  we follow the entropy method introduced  by  [GPV].   As explained in Section 1,  since  the   reference invariant state is not explicitly known,  we  compute the   entropy of the state of the process with respect to a 
 product measure with slowly varying profile $ \g (\cdot)$.
  We prove  in Lemma \eqv (dirichlet-l2)  that,  provided  $ \g (\cdot)$ is smooth enough, $ C^1$ suffices,     and takes the prescribed value $b (\cdot)$ at the boundary,    the rate to which the entropy  increases   is  of the order  of the volume, $ N^d$,  i.e  the same order of  the entropy and for finite   time  $T$   this  implies only a  modification of the constant multiplying  $ N^d$.
 
  We divide the proof of the hydrodynamic behavior 
in three steps: tightness of the measures $(Q_{\mu^N})_{N\geq 1}$, 
  energy estimates
and  identification of the  support of $Q^*$ as weak
solution of \equ(heq1) with fixed boundary conditions. We then refer to [KL], Chapter IV, that 
presents arguments, by now standard, to deduce the
hydrodynamic behavior of the empirical measures from the preceding 
results and the uniqueness of   the
weak solution of \equ(heq1).  We state without proving  the first two steps, tightness of the measures and   energy estimates. The proof of them can be  easily derived  from results already in the literature,    which we refer to. 
 
\medskip
\noindent{\bf \Proposition (lem1)} ({\bf Tightness})
{\it For almost any disorder configuration
$\a\in\S_D$,
the sequence $(Q_{\mu^N})_{N\geq 1}$  is tight and all 
its limit points
$Q^*$ are concentrated on absolutely continuous paths
$\pi(t,du)=\r(t,u)du$ whose density $\r$ is positive and bounded
above by $1$~:
$$
Q^*\Big\{ \pi\, :\, \pi(t,du)=\r(t,u)du\Big\}=1\; ,\quad
Q^*\Big\{ \pi\, :\, 0\le \r(t,u)\le 1\Big\}=1\; .  
\Eq(abscont)$$ } 

\medskip
Tightness  for non gradient  systems in contact with reservoirs is proven  in a way similar   to  the one for non gradient systems with periodic boundary conditions, see [KL], Chapter 7, Section 6.  The main difference relies on the fact that   for systems in contact with reservoirs   the invariant states are not product probability measures and some additional argument is required.
This can be proven as in    [LMS], Section 6.     
\medskip
In the next step we  
prove  
that for almost any disorder configuration
$\a\in\S_D$, every limit point $Q^*$ of the sequence 
$(Q_{\mu_N})_{N\geq 1}$  
is concentrated on paths whose densities $\r$ 
satisfy \equ(energy).  

\medskip
\noindent{\bf \Proposition (lemIB1)}
{\it For almost any disorder configuration
$\a\in\S_D$, every limit points $Q^*$ of the
sequence $(Q_{\mu^N})_{N\geq 1}$ is concentrated on the trajectories
that satisfies (IB1).
}
\medskip
The proof   can be done applying   arguments as in Proposition A.1.1. of  [KLO].
However the latter proof requires an application of Feynman-Kac formula, for which we have to replace our dynamic \eqv (gener1) (cf. [FM]).
\medskip 
 We then show that   $\P-$ a.s.  any  limit  point $Q^*$ 
 is supported on densities $\r$   satisfying 
\eqv (heq1) in the weak sense. 
 For  $\ell \in \N$, $x\in\L_N$,  with $ - N+ \ell \le  x_1 \le N- \ell $ denote by  $ \h^\ell (x)$  the average density of
$\h$ in a cube of width $2\ell+1$ centered at $x$   
$$ 
\h^\ell(x) = \frac 1{(2\ell+1)^d} \sum_{ y: |y-x|\le \ell}\h(y) . \Eq
(average) 
$$  
For a function $G$ on $\L$, $e\in\EE$, $\partial_{e}^N G$ denotes 
the discrete (space) derivative in the direction $e$
$$
\big(\partial_{e}^N G\big) (x/N)\;=\; N [ G((x+e)/N) - G(x/N)] \quad \hbox {with} \quad   x    \quad \hbox {and} \quad   x+e \in \L_N,  
\Eq (discrete) 
$$
and to short notation we denote by  $\partial_{k}^N G:=\partial_{e_k}^N G$ for $1\le k\le d$.

\medskip
\noindent{\bf \Proposition (lem2)}
{\it Assume that $D(\r)$ defined in \eqv (dif) can be
continuously extended in $[0,1]$. Then,  for almost any disorder 
configuration $\a\in\S_D$, any function $G$ in 
${\CC}_c^{1,2}([0,T]\times\iL)$ and any $\d>0$, we have 
$$  
\limsup_{c \to 0} \limsup_{a \to 0} \limsup_{N \to \infty} 
\Pr_{\mu^N} \left (
\left | \BB^{G,N}_{a,c} \right| \ge\d \right )=0, \Eq (ident-lim) 
$$
where
$$
\eqalign { 
&  \BB^{G,N}_{a,c}=
N^{-d}\sum_{x\in\L_N} G(T,x/N) \h_T(x) 
-  N^{-d} \sum_{x\in\L_N}G(0,x/N) \h_0(x)  
-  N^{-d}  \sum_{x\in\L_N} \int_0^T  \partial_s G(s,x/N)\h_s (x) ds \cr 
& \qquad  + \sum_{1\leq k,m\leq d} \int_0^T  ds N^{1-d}\sum_{x\in\L_N}  
\big(\partial_{k}^N G\big)(s, x/N)  \left \{ 
 D_{k,m} \left (\h_s^{[aN]}
(x)\right)  \right.  \cr    
&  \qquad \qquad \qquad \qquad \qquad \qquad \qquad \qquad \qquad \qquad 
 \left.  \times  \left\{ {(2c)}^{-1}
\Big[\h_s^{[aN]}(x+{c}Ne_m) 
-\h_s^{[aN]} (x-{c}N e_m) \Big]  
\right \}  \right\}. 
} \Eq (id-lim) 
$$
  }
\medskip \noindent 
The proof  is given  in Subsection 3.3. 
\medskip \noindent 
The last step states that   $\P-$ a.s.,  any limit points $Q^*$ of the
sequence $(Q_{\mu^N})_{N\geq 1}$ is  concentrated on the trajectories
with fixed density at the boundary and equal to $b(\cdot)$:

\medskip
\noindent{\bf \Proposition (lem2b)}
{\it $\P-$ a.s., 
any limit point  $Q^*$ of the
sequence $(Q_{\mu^N})_{N\geq 1}$ is concentrated on the trajectories
that satisfy (IB3).
}
 \medskip \noindent 
 The proof is given in  Subsection 3.4. 
 \medskip
\noindent{\sl 3.2. Basic estimates}
\medskip

\noindent{\bf \Lemma(lem-ergodic)} ({\bf Ergodic lemma}) 
{\it Let $V : \S_D\times \L \to \R$ a
bounded function, local with
respect to the first variable and continuous with respect to the
second variable, that is for any $\alpha\in \S_D$ the function $u\to
V (\alpha, u)$ is continuous and there exists an integer
$\ell \ge 1$ such that for all $u\in \L$ the support of
$V(\cdot,u) \subset \{-\ell,\cdots,\ell \}^d$. Then
$$
\lim_{N\to\infty} 
N^{-d} \sum_{x\in \L_N} \tau_xV(\alpha,x/N) 
=\int_\L \E \big[ V(\cdot ,u) \big] du\; 
\qquad  \P\ {\rm a.s.}.
\Eq(lem-1)$$
}

\smallskip
\noindent{\bf Proof.} 
We decompose the left hand side of the limit \eqv(lem-1) in two
parts
$$
\eqalign{ 
N^{-d} \sum_{x\in \L_N} \tau_x V (\alpha,x/N)  & 
= N^{-d} \sum_{x\in \L_N} 
  \Big( \tau_xV(\alpha,x/N) - \E \big[ V(\cdot ,x/N) \big] \Big)\cr
\ &+ N^{-d} \sum_{x\in \L_N} \E \big[ V(\cdot ,x/N) \big]
- \int_\L \E \big[ V(\cdot ,u) \big] du\; .
}
$$
By the stationary of $\P$ and the continuity of $u\to \E\big[ V(\cdot,u)\big]$, the second  term of
the the right hand side of the last equality converges to 0 as
$N\to\infty$. The first term converges to 0,
from Chebychef inequality and the classical method of
moments usually used in the proof of strong law of large numbers. 
\qed

\medskip

We start recalling the definition  of  relative  entropy, which is the main tool in the [GPV] approach. 
  Let  $
\nu_{\r(\cdot)}^{\a,N}$ be  the   product  measure defined in \eqv (inv1) and 
  $\mu$     a probability measure on $\SS_N$. Denote by
$H(\mu| \nu_{\r(\cdot)}^{\a,N} )$ the relative entropy of $\mu$ with
respect to $\nu_{\r(\cdot)}^{\a,N}$:
$$
H(\mu|\nu_{\r(\cdot)}^{\a,N}) \; =\; \sup_{f} \Big\{
\int f(\h) \mu (d\h) - \log \int e^{f(\h)} \nu_{\r(\cdot)}^{\a,N}  (d\h) \Big\}\; , 
$$
where the supremum is carried over all bounded
functions on $\SS_N$. Since $\nu_{\r(\cdot)}^{\a,N}$ gives a
positive probability to each configuration, $\mu$
is absolutely continuous with respect to $\nu_{\r(\cdot)}^{\a,N}$ and  we
have an explicit formula for the entropy: 
$$
H(\mu|\nu_{\r(\cdot)}^{\a,N} ) \; =\; \int  
\log  \Big\{ \frac{d\mu }{d\nu_{\r(\cdot)}^{\a,N}  } \Big\} \, d\mu 
\; . \Eq (ent1)
$$
Further, since  there is at most one particle per site, there exists a constant $C$, that depends only on  $\r(\cdot)$,  such that for all $\a$
$$
H(\mu| \nu_{\r(\cdot)}^{\a,N} ) \;\leq\; C  N^d
\Eq (entbound)$$
for all probability measures $\mu$ on $\SS_N$  (cf. comments
following Remark V.5.6 in [KL]).

It is well known that  one of the main step in the derivation of hydrodynamic limit for the empirical density is a super exponential estimate which allows the replacement of local functions by functionals of the empirical density. One needs to estimate expression such as $<Z,f>_{\mu^N}$ in terms of Dirichlet form $ <- 
\LL_N   \sqrt{f(\h)},  \sqrt{f(\h)} >_{\mu^N}$, where $Z$ is a local function and $<\cdot,\cdot>_{\mu^N}$ represents a scalar product with respect to some state $ \mu^N$.  Since in the context of boundary driven process the invariant state is not explicitly known and   we fix as reference measure some product measure $\nu$, see Lemma \eqv (dirichlet),   there are no reasons for  $ <- 
\LL_N   \sqrt{f(\h)},  \sqrt{f(\h)} >_{\nu}$ to be positive.  Next lemma shows that this expression is almost positive.
Let  $\DD_N^{0} ( \cdot, \nu)$, $\DD_N^{b} ( \cdot, \nu)$ be functionals
from    $h \in L^2 ( \nu)$ to $\R^+$: 
$$
\eqalign{
\DD_N^{0}\big(h, \nu\big) & =
\frac12 \sum_{e\in\EE} \sum_{x,x+e\in\L_N } \int
     C_N (x,x+e;\eta)
  \left( {h}(\eta^{x,x+e}) -{h}(\eta) \right)^2 d \nu (\eta) \, ,\cr
\DD_N^{b} \big(h,\nu) & = 
\frac12 \sum_{x\in\Gamma_N} \int  C^{b}(x/N,\eta)
  \left( {h}(\eta^{x}) -{h}(\eta) \right)^2 d \nu (\eta) \, .
}
\Eq (dir-form2)
$$
 
\medskip
\noindent {\bf \Lemma (dirichlet)} 
{\it Let $\gamma:\L\to (0,1)$ be a smooth function such that
$\g{\big\vert_\Gamma} \; =\; b(\cdot)$.
For any
$\a\in \S_D$ and $a >0$
there exists a positive 
constant $C_0\equiv C_0(A, \|\nabla \g \|_\infty )$  so  that for any 
$f\in L^2\big( \nu_{\g(\cdot)}^{\a,N} \big)$,       
 
$$
 \int_{\SS_N} f(\h) \LL^0_{N}
f(\h)  d  \nu_{\g(\cdot)}^{\a,N}(\h) 
\le -\big(1-\frac{1}{2a}\big) {\DD}_N^{0} \big({f},\nu_{\g(\cdot)}^{\a,N} \big)  
+C_0 N^{d-2}(a+1)  \| f\|^2_{L^2(\nu_{\g(\cdot)}^{\a,N}) },   \Eq (bound3)  $$ 
$$ \int_{\SS_N} f(\h) \LL^b_{N}
f(\h)  d  \nu_{\g(\cdot)}^{\a,N}(\h)  
= -  {\DD}_N^{b} \big({f},\nu_{\g(\cdot)}^{\a,N} \big) \; . 
 \Eq (bound4) 
$$
}
\medskip
\noindent{\bf Proof.}
By  \equ(dir-form2) ,
$$
\eqalign{
&\int_{\SS_N} f(\h)\LL^0_{N}
f(\h)  d  \nu_{\g(\cdot)}^{\a,N}(\h)
=- {\DD}_N^{0}({f},\nu_{\g(\cdot)}^{\a,N}) \cr  
&\qquad\qquad
 + \frac12 \sum_{e\in\EE} \sum_{x,x+e\in\L_N} \int
C_N(x,x+e;\eta) \big(\nabla_{x,x+e}f \big) (\eta) f(\eta^{x,x+e})
R_1(x,x+e;\eta)
d  \nu_{\g(\cdot)}^{\a,N}(\h) \; , 
}
$$
where 
$$
R_1(x,x+e;\eta) = \big(\nabla_{x,x+e}\eta(x) \big)
\big( e^{(N^{-1} \partial_e^N \l_0(\g(x/N))  )} -1  \big)\; . 
$$
By the elementary inequality $2 u v\leq a  u^2 +a^{-1} v^2$ 
which holds for any    $a>0$, for any $x,x+e\in \L_N$
$$\eqalign{ &
\int C_N(x,x+e;\eta)   (\nabla_{x,x+e}{f}) f(\eta^{x,x+e}) R_1(x,x+e,\eta ) d
\nu_{\g(\cdot)}^{\a,N}(\h)  \cr 
&\quad
\leq \frac{1}{2 a}\int C_N(x,x+e;\eta)   (\nabla_{x,x+e} {f})^2 d  \nu_{\g(\cdot)}^{\a,N}(\h)
+\frac{a}{2}\int   C_N(x,x+e;\eta) f(\eta^{x,x+e})^2 (R_1(x,x+e))^2  d
  \nu_{\g(\cdot)}^{\a,N}(\h)\; .  \cr 
}
$$
To conclude the proof it remains to use Taylor expansion and an
integration by part   in the second term of the right hand side of the last inequality.
On the other hand, since $\g{\big\vert_\Gamma} = b(\cdot)$ the
measure $\nu_{\g(\cdot)}^{\a,N}$ is reversible with respect to
$\LL^b_{N}$.  A simple computation shows that
$$
\int_{\SS_N} f(\h) \LL^b_{N} f(\h)  d  \nu_{\g(\cdot)}^{\a,N}(\h)
=- {\DD}_N^{b}({f} ,\nu_{\g(\cdot)}^{\a,N} ) \;  .$$
\qed

\medskip
\noindent {\bf \Lemma (dirichlet-l1)} 
{\it Let $\rho,\rho_0:\L\to (0,1)$ be two smooth functions.
There exists a positive 
constant $C_0'\equiv C_0'(A, \|\nabla \rho_0\|_\infty ,\|\nabla \rho\|_\infty )$ such that for any probability  measure $\mu^N$
on $\SS_N$ and for any $\a\in \Sigma_D$,     
$$
 \DD_N^{0} \Big(\sqrt{ \frac{d\mu^N}{d\nu_{\rho(\cdot)}^{\a,N}} },
\nu_{\rho(\cdot)}^{\a,N}  \Big)\; \le \;
2\; \DD_N^{0} \Big(\sqrt{ \frac{d\mu^N}{d\nu_{\rho_0(\cdot)}^{\a,N}} },
\nu_{\rho_0(\cdot)}^{\a,N}  \Big)
  \; +\; C_0' N^{d-2}\; .
 \Eq (bound-l1) 
$$
}
\medskip

\noindent{\bf Proof.}
Denote by $f(\eta)=\frac{d\mu^N}{d\nu_{\rho(\cdot)}^{\a,N}}(\eta) $
and $h(\eta)=\frac{d\mu^N}{d\nu_{\rho_0(\cdot)}^{\a,N}}(\eta)$. Since  $f(\eta)=  h(\eta)  \frac { {d\nu_{\rho_0(\cdot)}^{\a,N}}(\eta) }  {d\nu_{\rho(\cdot)}^{\a,N}(\eta)} $ we obtain  for
$e\in\EE$ and $x,x+e\in\L_N$ the following  
$$\eqalign{
&\int_{\SS_N} C_N(x,x+e;\eta) \Big[ \nabla_{x,x+e}\sqrt{f}(\eta) \Big]^2
d\nu_{\rho (\cdot)}^{\a,N}(\eta) \cr
\ & \ \ 
=\int_{\SS_N} C_N(x,x+e;\eta) \Big[ \sqrt{h}(\eta^{x,x+e})
R_2(x,x+e;\eta) +\nabla_{x,x+e} \sqrt{h}(\eta) \Big]^2
d\nu_{\rho_0(\cdot)}^{\a,N}(\eta)  \cr
&\ \ \ 
 \le 2 \int_{\SS_N} C_N(x,x+e;\eta) \Big[ \nabla_{x,x+e}\sqrt{h}(\eta) \Big]^2
d\nu_{\rho_0 (\cdot)}^{\a,N}(\eta) \cr
\ &\quad
+2\int_{\SS_N} C_N(x,x+e;\eta) h(\eta^{x,x+e})
\big[R_2(x,x+e;\eta)\big]^2 
d\nu_{\rho_0(\cdot)}^{\a,N}(\eta) \; ,
}
$$
where 
$$
R_2(x,x+e;\eta) = 
\exp \big\{ (1/2)N^{-1} \partial_e^N [ \l_0(\rho(x/N))-\l_0(\rho_0(x/N))  ]  \nabla_{x,x+e}\eta(x) \big\} -1 \; . 
$$
We conclude the proof using Taylor expansion and integration by
parts.
\qed
 
\medskip
Denote by $S_t^N$ the semigroup
associated to the generator $N^2 \LL_N$.
Given  a   probability measures $\mu^N$ on $\SS_N$ denote  by $\mu^N (t)$ the state of
the process at time $t$~: $\mu^N (t) =\mu^N S_t^N$.

Recall that $\gamma\colon\Lambda\to (0,1)$ is a smooth profile equal to
$b$ at the boundary of $\Lambda$. Let $h_t^N$ be the density
of $\mu^{N} (t)$ with respect to $\nu_{\gamma(\cdot)}^{\a,N}$.  
Let $  \LL_{\g,N}^*$ be the adjoint of $\LL_N$ in
$L^2(\nu_{\gamma(\cdot)}^{\a,N})$. It  is easy to
check that
$$
\partial_t h_t^N \; =\; N^2 \LL_{\g,N}^* h_t^N \;.\Eq(Chapman)
$$
 Notice that $\LL_{\gamma, N}^*$ is not a
generator because $\nu_{\gamma(\cdot)}^{\a,N}$ is not an invariant measure
for the Markov process with generator $\LL_N$. We denote by $H_N(t)$ the
entropy of $\mu^N(t)$ with respect to $\nu_{\gamma(\cdot)}^{\a,N}$, see \eqv (ent1), 
$$ H_N(t):= H(\mu^N (t)|\nu_{\gamma(\cdot)}^{\a,N} ). \Eq (ent2) $$

\medskip
\noindent {\bf \Lemma (dirichlet-l2)} 
{\it    There exists positive constant $C=C(\|\nabla \gamma\|_\infty)$  such that  for any $a>0$ and for any  $\a\in \S_D$
$$
 \partial_t  H_N(t) \; \le \; -  2(1-a)  N^2 
\DD_N^{0}  (\sqrt {h_t^{N}}, \nu_{\gamma(\cdot)}^{\a,N} )
 -  2 N^2 
\DD_N^{b}  (\sqrt {h_t^{N}}, \nu_{\gamma(\cdot)}^{\a,N} )\; +\;
\frac{C}{a} N^d\; ,
$$
}
\medskip

\noindent{\bf Proof.}  By 
\equ(Chapman) and the explicit formula for the entropy we have that
$$
\partial_t H_N(t) =N^2 \int_{\SS_N} h_t^N \LL_N \log \big( h_t^N  \big) d \nu_{\gamma(\cdot)}^{\a,N}  \; .
$$
Using the basic inequality
$
a\big(\log b -\log a  \big)\le -\big( \sqrt{a}-\sqrt{b}\big)^2+\big(b-a\big)
$
for positive $a$ and $b$, we obtain
$$\eqalign{
\partial_t H_N(t) & \; \le\;
-2N^2 \DD_N^{0} \big(\sqrt {h_t^{N}},\nu_{\gamma(\cdot)}^{\a,N}  \big) 
-2N^2 \DD_N^{b} \big(\sqrt {h_t^{N}},\nu_{\gamma(\cdot)}^{\a,N}   \big) \cr
\ &\ \ \ 
+N^2\int_{\SS_N}  \LL^0_{N} h_t^N  d \nu_{\gamma(\cdot)}^{\a,N}  
+ N^2\int_{\SS_N} \LL^b_{N} h_t^N  d \nu_{\gamma(\cdot)}^{\a,N}  \; .
}
\Eq(entr1)
$$
Since $\g  (u)= b(u)$  for $u \in \G$,  $\nu_{\gamma(\cdot)}^{\a,N} $
is reversible with respect to $\LL^b_{N}$.  This implies that 
$$\int_{\SS_N}\LL^b_{N}
h_t^N  d \nu_{\gamma(\cdot)}^{\a,N} =0.
$$
We shall now obtain a bound for 
$\int_{\SS_N}  \LL^0_{N} h_t^N  d \nu_{\gamma(\cdot)}^{\a,N}  $ in terms of
$\DD^0_N$.
Denote by $R:\R\to \R$ the function defined by
$
R(u)=e^u-1-u \; .
$
A standard computation shows that 
$$
\eqalign{
&N^2\int_{\SS_N} \LL^0_{N} h_t^N  d \nu_{\gamma(\cdot)}^{\a,N}  \cr
&\qquad 
=N^2\sum_{e\in\EE}\sum_{x,x+e\in \L_N} \int C_N(x,x+e;\eta) h_t^N(\eta) 
R\big(N^{-1} \partial_e^N \l_0(\g (x/N))\nabla_{x,x+e}\eta(x)  \big) 
 d \nu_{\gamma(\cdot)}^{\a,N}  (\eta)\cr
 &\qquad \quad + N\sum_{e\in\EE}\sum_{x,x+e\in \L_N} (\partial_e^N \l_0(\g
(x/N))  \int   W_{x,x+e}(\eta) h_t^N(\eta) d
\nu_{\gamma(\cdot)}^{\a,N} (\eta)\; ,
}\Eq(entr2)$$ 
where  
$W_{x,x+e}(\h)$ is the current  over the bond $(x,x+e)$ :
$$
W_{x,x+e}(\h)\equiv
 C_N(x,x+e;\h)\big[\h(x)-\h(x+e)\big]
\; . 
\Eq(current)$$
We will often omit to write the  dependence of
$W_{x,x+e} (\eta) $
on $N$ and $\h$. 
By Taylor expansion and the elementary
inequality $|R(u)|\le \frac{u^2}{2}e^{|u|}$, we
obtain using the fact that $\g$ is  smooth and $h_t^N$ is a probability
density with respect to $\nu_{\gamma(\cdot)}^{\a,N} $, that the first term
of the right hand side of the \equ(entr2) is bounded by $C\; N^d$
for some positive constant $C$.  On the other hand   integrating by
part,  applying  the same computations as in  Lemma 5.1  of [LMS],  we obtain  that there exists a   constant  $C_0= C(\|\nabla \gamma\|_\infty)$  so that  for any $a>0$ 
$$
\int W_{x,x+e} h_t^N d \nu_{\gamma(\cdot)}^{\a,N}  \; \le\;
\frac{1}{a} \int C_N(x,x+e;\eta) \left( \nabla_{x,x+e}\sqrt{h_t^N} \right)^2
d\nu_{\gamma(\cdot)}^{\a,N}  
+C_0 \big\{a+  N^{-1}\big\}
$$
for  $x,x+e\in \L_N$. \qed

\medskip

For $z \in \L_N$,  $M\in \N$   denote  by $\Lambda_M (z) $ the intersection of  a cube centered at   $z \in \L_N$ of edge $2M+1$ with $ \L_N$, i.e 
$$  \Lambda_M (z):=   {\{z + \Lambda_M\} }\cap  \L_N. \Eq (cube.1)$$
For  probability measure
$\nu^N$ on $\SS_N$, denote by $\DD_{M,z}^{0}
(\cdot\, , \nu^N)$ the Dirichlet form corresponding to
jumps  in  $\Lambda_M (z)$:
$$
\DD_{M,z}^{0} (f, \nu^N) \; =\; \frac 12  \sum_{x,  x+e \in \Lambda_M (z)} 
\int C_N(x,x+e;\eta)  (\nabla_{x,x+e}  f (\eta) )^2
d\nu^N (\eta) \; .  \Eq (D.1)
$$
 Similarly, for $z\in \Gamma_N$
define $\DD_{M,z}^{b} (\cdot\, , \nu^N)$   the Dirichlet form
corresponding to creation and destruction of particles at sites in  $\Gamma_N$  which are at distance less than $M$ from $z$~:
$$
\DD_{M,z}^{b} (f, \nu^N) \; =\; \frac 12 \sum_{\scriptstyle x\in
  \Gamma_N \cap  \Lambda_M (z) } 
\int C^{b} (x/N,\eta) \big( f  (\eta^x) -f  (\eta)\big)^2
d\nu^N (\eta)  \; .  \Eq (D.2)
$$
Fix any   $z \in \Gamma_N$  denote by $f_{t}^{z,N}$   the Radon-Nikodym derivative of $\mu^N
(t)$ with respect to $\nu_{ b(z/N)}^{\a,N} $, the random  Bernoulli measure on $\SS_N$   with constant parameter equal to  $b(\frac z N)$.  Recall that  we denoted by $h_t^N$   the Radon-Nikodym derivative of $\mu^N
(t)$ with respect to $\nu_{\g(\cdot)}^{\a,N} $ and that  $b(\frac z N)= \g(\frac z N)$ for $z \in \G$.  We have the following result.  
\medskip
\noindent {\bf \Lemma (dirichlet-l2b)} 
{\it  Take  $ M \in \N$, $M< N$.
There exists a positive constant  $C_0= C(\|\nabla\gamma\|_\infty )$ depending only on  $\gamma(\cdot
)$ such that for any $z\in \Gamma_N$  
$$
\eqalign{
 \DD_{M,z}^{0} \big(\sqrt {f_t^{z,N}}, \nu_{b(z/N)}^{\a,N} \big) &\; 
 \le \;  2  \DD_{M,z}^{0} \big(\sqrt {h_t^N} ,
      \nu_{\gamma(\cdot)}^{\a,N} \big) + C_0 \frac{M^d}{N^2} \; ,\cr
 \DD_{M,z}^{b} \big(\sqrt {f_t^{z,N}}, \nu_{b(z/N)}^{\a,N}  \big) &\; \le \;  
                     2 \DD_{M,z}^{b} \big(\sqrt {h_t^N} ,
                      \nu_{\gamma(\cdot)}^{\a,N} \big) + C_0
      \frac{M^{d+1}}{N^2}\; . \cr                     
}
$$
}
\medskip
\noindent  The proof is  similar to the proof
of Lemma \equ(dirichlet-l1).  
\medskip
\noindent{\sl 3.3. Proof of Proposition \equ(lem2)}
\medskip 
We prove in this section Proposition \equ(lem2). 
Let  $Q^*$  be  a limit point of the sequence 
$(Q_{\mu^N})_{N\ge 1}$ and assume, without loss of
generality, that $\P-$ a.s.,  $Q_{\mu^N}$ converges to $Q^*$.
Fix a function $G$ in ${\CC}_c^{1,2}([0,T]\times\iL)$. For  
$\a\in \O_D$  consider the   $\Pr_{\mu^N}$ martingales  with
respect to the 
natural filtration associated with $(\h_t)_{t\in[0,T]}$,
$M_t^G\equiv M_t^{G,N,\a}$ and  $\NN_t^G \equiv \NN_t^{G,N,\a}$, 
$t\in[0,T]$, defined by 
$$ 
\eqalign{
M_t^G &\; =\; <\pi_t^N,G_t>-<\pi_0^N,G_0> -\int_0^t \big (  <\pi_s^N,
\partial_s G_s>+N^2
\LL_N <\pi_s^N, G_s> \big)  \, d s \; , \cr
 \NN_t^G &\; =\; \left(M_t^G \right)^2 
 -\; \int_0^t \left\{ N^2 \LL_N^{\a}
\big(<\pi_s^N,G_s>\big)^2 - 2 <\pi_s^N,G_s> N^2 \LL_N
<\pi_s^N,G_s>\right\} d s \; .
} \Eq  (MMa1)
$$ 
A computation of the integral term of $\NN_t^G$ shows
that the expectation of the quadratic variation of 
$M_t^G$ vanishes as $N\uparrow 0$. Therefore, by Doob's
inequality, for every $\d >0$, $\P=1$, 
$$ 
\lim_{N\rar \infty} \Pr_{\mu_N} \Big[ \sup_{0\le t\le T}
|M_t^G|>\d \Big] \; =\;0 \; . 
\Eq(bv)
$$ 
Thanks to \eqv (C.3) and since for any $s\in[0,T]$ the  function
$G_s$ has  compact support in $\iL$, a summation by
parts permits to rewrite the integral term of $M_t^G$ as
$$ 
\int_0^t <\pi_s^N,\partial_s G_s> d s \; +\; \int_0^t  \, 
\Bigl\{ N^{1-d}\sum_{k=1}^{d} \sum_{x\in \L_N} \big(\partial_{k}^N
G_s \big) (x/N)  W_{x,x+e_k}(\h_s) \Bigr\}ds , 
\Eq (E.11) $$ 
where the current $W_{x,x+e_k}$ is defined in \equ(current).
To localize the dynamics 
define
for any $0<r<1$  
$$
\eqalign { & \L_{r} = [-r,r]\times \T^{d-1} ,\qquad \L_{rN} =\{(x_1,\cdots,x_d)\in\L_N\ :\ \ -rN \le x_1\le rN \},  \cr &
 \G_{rN} =\{x \in\L_{rN} \ :\ \  x_1= \pm rN \}.  } 
\Eq(subsets)
$$
Set, for $0<a <c<1$, $k=1,\ldots,d$,
$$
\V_{k}^{N,c,a} (\h,\a)\; = \; 
N W_{0,e_k} 
 + \sum_{m=1}^{d} D_{k,m} \left(\h^{[aN]}(0)\right)
 \left\{ {(2c)}^{-1}
\Big[ \h^{[aN]}({c} N e_m) -\h^{[aN]} (-{c}Ne_m) 
\Big] \right\}.
 \Eq (E.21) 
$$
Next theorem is the main step in the proof of Proposition \equ(lem2).

\vskip0.5cm
  \noindent{\bf \Theorem (thng2)}
{\it  Assume that $D(\cdot)$ defined in \eqv (dif) can be 
continuously extended in $[0,1]$. Then, $\P=1$,   for any $G\in {\CC}_c^{1,2}([0,T]\times\iL)$,
$$ 
\limsup_{c \rar 0}
\limsup_{a \rar 0}\limsup_{N \rar \infty}
\Es_{\mu_N} \Big[ \, \Big| N^{-d} \int_0^T 
 \sum_{x\in \L_N}G_s(x/N)\tau_x \V_{k}^{N,c,a} (\h_s,\a)\, d s
\Big|\, \Big] =0 \Eq (E.24)
$$ 
for $k=1,\ldots,d$. 
}

\medskip
\noindent{\bf Proof.}    Let $0< \theta <1$ such that for any $t\in[0,T]$, the
support of the function $G_t$ is a subset of
$\Lambda_{(1-2\theta)}$.  
Fix a smooth function $\gamma_\theta\colon \Lambda\to (0,1)$ which
coincides with $b$ at  the boundary of $\Lambda$ and  constant inside
$\L_{(1-\theta)}$. 
Denote by $Z_k^{N,c,a} (G,\eta)$ the
quantity
$$
Z_k^{N,c,a}\big(G,\eta\big) \; =\; N^{-d}
\sum_{x\in \L_N}G(x/N)\tau_x  \V_{k}^{N,c,a} (\h_s,\a)\; .
$$
Since the entropy  of $\mu^N$ with respect to $\nu_{\gamma_\theta(\cdot)}^{\a,N}$ 
is bounded by $C_\theta |\L_N|$ for some finite constant $C_\theta$, by the entropy 
inequality, the left hand side of \equ(E.24) is bounded above by
$$
\frac{C_\theta}{B} \; + \; \frac{1}{B N^d}\log 
\Es_{\nu_{\gamma_\theta (\cdot )}^{\a,N} } 
\Big[ \exp \Big\{ B N^d \Big| \int_0^T Z_k^{N,a,c}
\big( G_s,\eta_s\big)d s\Big| \Big\} \Big] \Eq (in1)
$$
for any positive $B$.
Since $e^{|x|}\le e^{x}+e^{-x}$ and $\limsup N^{-d} \log \{a_N + b_N\}
\le  \max \{ \limsup N^{-d} \log a_N$ , $\limsup N^{-d} \log b_N \}$,
we may remove the absolute value in the second term of \eqv (in1), provided our estimate 
remains in force if we replace $G$ by $-G$. By the Feynman-Kac 
formula,
$$
  \frac{1}{B N^d}\log 
\Es_{ \nu_{\gamma_\theta (\cdot )}^{\a,N}} 
\Big[ \exp \Big\{ B N^d  \int_0^T Z_k^{N,a,c}
\big( G_s,\eta_s\big)d s  \Big\} \Big]  \le \frac{1}{ B N^d}\int_0^T \lambda_{N,c,a}(G_s) \, ds\; ,
$$
where $\lambda_{N,c,a}(G_s)$ is the largest eigenvalue of the 
$N^2\{\LL_N^{sym} + B
Z_k^{N,c,a}( G_s,\eta)\}$  where  $  \LL_N^{sym}: = \frac 1 2 (  \LL_N +  \LL^*_{\g_\th,N} )$ and  $  \LL^*_{\g_\th,N} $ is the  adjoint of $ \LL_N $ in $L^2(\nu_{\gamma_\theta (\cdot )}^{\a,N})$. By the variational formula for the
largest eigenvalue, for  $s \in [0,T]$, we have that 
$$
 \frac 1 {B N^d} 
\lambda_{N,c,a}(G_s) = \sup_{f} \Big\{\int Z_k^{N,c,a}\big(G_s,\eta\big) f(\eta)
\nu_{\gamma_\theta(\cdot)}^{\a,N}(d\eta) \; +\; \frac{N^{2-d}}{B} < \LL_N \sqrt{f} , 
\sqrt{f} >_{\gamma_\theta (\cdot)} \Big \}\; .
$$
In this formula the supremum is carried over all densities $f$
with respect to $\nu_{\gamma_\theta (\cdot )}^{\a,N}$ and notice that we used $ < \LL_N \sqrt{f} , 
\sqrt{f} >_{\gamma_\theta (\cdot)}=  <  \LL_N^{sym} \sqrt{f} , 
\sqrt{f} >_{\gamma_\theta (\cdot)} $. 
Since $\gamma_\theta (\cdot)$ coincides with $b(\cdot)$  on $\G$, 
$\LL^b_{N}$ is reversible with respect to $\gamma_\theta (\cdot)$, so that 
 $< \LL^b_{N} \sqrt{f} , 
\sqrt{f} >_{\gamma_\theta(\cdot)}$ is negative. We then apply  simply \eqv (bound3)
of  Lemma \equ(dirichlet) with $a=1$ to estimate $< \LL_N \sqrt{f} , 
\sqrt{f} >_{\gamma_\theta(\cdot)}$ by $- (1/2) \DD_N^{0}(\sqrt{f},
\nu_{\gamma_\theta (\cdot )}^{\a,N}) + C_\theta'N^{d-2}$ for some
constant $C_\theta'$. In particular, to prove the theorem, we just need to show that
$$
\limsup_{c \rightarrow 0}\limsup_{a \rightarrow 0}
\limsup_{N \rightarrow \infty}
\int_0^T ds \sup_{f} \Big\{ \int Z_k^{N,c,a}\big( G_s,\eta\big) f(\eta)
\nu_{\gamma_\theta(\cdot)}^{\a,N}(d\eta) - 
\frac 1 B  N^{2-d} \DD_N^{0} (\sqrt{f}, \nu_{\gamma_\theta}^{\a,N})
\Big \} \;=\; 0  
$$
for every $B>0$ and then let $B\uparrow\infty$.
Notice that for $N$ large enough and $a,c$ small enough, the function 
$Z_k^{N,c,a}(G_s,\eta)$ depends on the configuration $\eta$ only through
the variables $\{ \eta(x),\; x\in \L_{(1-\theta)N}\}$. Since
$\gamma_\theta (\cdot)$ is constant, say equal to $\gamma_0$  in
$\L_{(1-\theta)}$, we may replace $\nu_{\gamma_\theta (\cdot)}^{\a,N}$
in the previous formula by $\nu_{\gamma_0}^{\a,N}$.  The $\nu_{\gamma_0}^{\a,N}$ is
reversible for $\LL^0_{N}$ and therefore $\DD_N^{0}(\cdot\, ,
\nu_{\gamma_0}^{\a,N})$ is the Dirichlet form associated to the generator
$\LL^0_{N}$. Since the Dirichlet form is convex, it remains to show that
$$
\limsup_{c \rightarrow 0}\limsup_{a \rightarrow 0}
\limsup_{N \rightarrow \infty}
\int_0^T ds \sup_{f} \Big\{ \int Z_k^{N,c,a}\big( G_s,\eta\big) f(\eta)
\nu_{\gamma_0(\cdot)}^{\a,N}(d\eta) - 
\frac 1 B  N^{2-d}\DD_N^{0} (\sqrt{f}, \nu_{\gamma_0 }^{\a,N})
\Big \} \;=\; 0  
$$
for every $B>0$. This result  has been proved in [FM], Theorem 3.2.  
\qed
 \medskip
\noindent {\bf Proof of Proposition  \eqv  (lem2):}     By   \eqv (MMa1), \eqv (E.11) and \eqv (E.21), applying Theorem \eqv  (thng2) we obtain   \eqv (ident-lim). \qed 

 \medskip
\noindent{\sl 3.4. Proof of Proposition \equ(lem2b)}
\medskip 

For $a>0$,  $u \in \Lambda$   denote    
$$
\iota_a (u)=  \frac{1}{ |\big[-a,a\big]^d \cap \Lambda |} \1_{\big\{\big[-a,a\big]^d \cap \Lambda \big\}}(u);      \Eq (identite)$$
and  for $A\subset \L$   define
 the sets $A^\pm$  as 
$$
A^+ =\{(u_1,\ldots,u_d)\in A\ \ :\ \ u_1>0 \},\quad A^-
=\{(u_1,\ldots,u_d)\in A\ \ :\ \ u_1<0 \}\; . \Eq  (set1)
$$
We define similarly $A^+_N$ and $A^-_N$ when  $ A_N \subset \L_N$. 
 Let $G (\cdot, \cdot) \in C^{1,2} \left ([0,T] \times \L \right ) $,
 $\mu \in D([0,T],\MM_1(\L))$ and 
for $0<a<c<1$, define the following functional 
$$\eqalign{
\hat F_{a,c}^{G}\big(\mu (\cdot, \cdot) \big)& =
 \int_0^T  ds   \int_{\L(1-c)}  du \Big\{  G_s (u)\;   {(2c)}^{-1}
\Big[ \big(\mu_s \star \iota_a\big)(u+{c}e_1) -\big(\mu_s \star \iota_a\big) (u-{c}e_1) 
\Big]  \Big\}\cr
\ &\quad
+\int_0^T  ds   \int_{\L} du \partial_{e_1}G_s(u) \big(\mu_s \star \iota_a\big)(u)
-\int_0^T ds\Big\{  \int_\Gamma b(u) {\hbox{\bf
n}}_1(u)  G_s (u) 
 \hbox{d} \hbox{S} \Big\}\; , 
}
 \Eq (funct1) $$
where $G_s(u)\equiv  G(s,u)$,  {\bf n}=$(\hbox{\bf n}_1,\ldots ,\hbox{\bf
  n}_d)$ is the outward unit normal vector to the boundary surface
$\Gamma$ and $\hbox{d} \hbox{S}$ is  the  surface element  of $\Gamma$.
The proof of Proposition \equ(lem2b) follows from the next lemma.

\medskip
\noindent {\bf \Lemma (bord)} 
{\it 
For $G (\cdot, \cdot) \in C^{1,2} \left ([0,T] \times
\L \right ) $, $\P$ a.s. we have 
$$ \limsup_{c \to 0} \limsup_{a \to 0}\limsup_{N \to \infty} 
\Es^{Q_{\mu^N}}\left [ \left|
\hat F_{a,c}^{G}\big(\mu^N (\cdot, \cdot) \big)\right| \right]  = 0. $$ 
}
\medskip
\noindent{\bf Proof.}  To short notation, denote $ f_s(u):= (\mu_s \star \iota_a\big)(u)$. Taylor expanding we have that 
$$ \eqalign { &  \int_{\L(1-c)}  du \Big\{  G_s (u)\;   {(2c)}^{-1}
\Big[ f_s(u+{c}e_1) -f_s (u-{c}e_1) 
\Big]  \Big\}  \cr & =
  \frac 1 {2c} \int_{(\Lambda \setminus
\L_{(1- 2c)})^+}    G_s(u-c e_1)  f_s(u) du-  
\frac 1 {2c} \int_{(\Lambda \setminus
\L_{(1-2c)})^-}     G_s(u+c e_1)  f_s (u) du 
\cr &-\int_{\L (1-2c)}     \partial_{e_1}  G_s (u) f_s(u)  du  +  c  \int_{\L (1-2c)}    R(G,c,s,u) f_s(u)  du.
 }  \Eq (t1) $$
where  $|R(G,c,s,u)| \le  \sup_{u \in \Lambda} \sup_{s \in[0,T] }|\partial^2_{e_1}G_s (\cdot)| $.
Since $f_s(u)\le 1$ uniformly in $s$ and $u$ 
$$\left | \int_{\L (1-c)}    R(G,c,s,u)  f_s(u)  du\right | \le  2 \sup_{u \in \Lambda} \sup_{s \in[0,T] } |\partial^2_{e_1}G_s (u)|,
 \Eq (bound8) $$
 and 
 $$  \left | \int_{\L (1-2c)}     \partial_{e_1}  G_s (u) f_s(u)  du  -\int_{\L}     \partial_{e_1}  G_s (u) f_s(u)  du \right |  \le  2c  \sup_{u \in \Lambda} \sup_{s \in[0,T] } |\partial_{e_1}G_s (u)|.  $$
 Taking in account \eqv (t1), \eqv (identite)  and \eqv (bound8)
the lemma is then  proven once we show that $ \P=1 $ the following holds
$$\eqalign{
&\limsup_{c \to 0} \limsup_{a \to 0}\limsup_{N \to \infty} 
\Es_{\mu^N}\Big[ 
\Big|\int_0^Tds \Big\{ \frac{1}{2c N^d} \sum_{x\in (\L_{(1-a)N}\setminus
   \L_{(1-a-2c)N})^\pm} G_s
 (\frac x N) \eta_s^{aN} (x) \cr
&\qquad\qquad\qquad\qquad\qquad\qquad\qquad\qquad\qquad
- \frac{1}{N^{d-1}}
 \sum_{x\in \Gamma_N^\pm} b(\frac x N) G_s (\frac x N) \Big\}\Big|\Big]=0\; ,
 }\Eq(bord2)
$$
where for $0<\ve< 1$, $\L_{\ve N}$ and   $(\L_{\ve N})^+$ are defined in  \equ(subsets) and  below  \eqv  (set1).
 By adding and subtracting the same quantity 
in the expectation of \equ(bord2), it is easy to see that the limit \equ(bord2) follows  once the 
next two lemmas are proven. \qed
\medskip
\noindent {\bf \Lemma (bord-l2)} 
{\it 
For $G (\cdot, \cdot) \in C^{1,2} \left ([0,T] \times
\L \right ) $, $\P$ a.s. we have 
$$\eqalign{
&\lim_{\ell\to\infty} \limsup_{c \to 0} \limsup_{a \to 0}\limsup_{N \to \infty} 
\Es_{\mu^N}\Big[ 
\Big|\int_0^Tds \Big\{ \frac{1}{2c N^d} \sum_{x\in (\L_{(1-a)N}\setminus
    \L_{(1-a-2c)N})^\pm} G_s (x/N) \eta_s^{aN} (x) \cr
&\qquad\qquad\qquad\qquad\qquad\qquad\qquad\qquad\qquad\qquad
- \frac{1}{N^{d-1}}\sum_{x\in\G_{ (1- \frac {\ell} N)N}^\pm} G_s(x/N) \eta_s^\ell(x) 
 \Big\}\Big|\Big]=0\; .
} \Eq(bord3)$$ 
}
\medskip

\noindent {\bf \Lemma (bord-l3)} 
{\it 
For $G (\cdot, \cdot) \in C^{1,2} \left ([0,T] \times
\L \right ) $, $\P$ a.s. we have 
$$\eqalign{
&\lim_{\ell\to\infty} 
\limsup_{N \to \infty} 
\Es_{\mu^N} \Big[ 
\Big|\int_0^Tds \Big\{ 
 \frac{1}{N^{d-1}}\sum_{x\in\G_{ (1- \frac {\ell} N)N}^\pm} G_s(x/N) \eta_s^\ell(x)  \cr
&\qquad\qquad\qquad\qquad\qquad\qquad\qquad\qquad\qquad\qquad
- \frac{1}{N^{d-1} }
 \sum_{x\in \Gamma_N^\pm} b(x/N) G_s (x/N) \Big\}\Big|\Big]=0\; .
}\Eq(bord4)
$$ 
}
\medskip

\noindent {\bf Proof of Lemma \equ(bord-l2). }
The summation in \equ(bord3) contains two similar terms. We consider
the one corresponding to the summation of the right hand side of
$\L_N$  (i.e. the one with signe +). 
By Taylor expansion applied to the function $G$, the expectation in
the statement of the lemma is  bounded above by
$$\Es_{ \mu^N} \Big[ 
\Big|\int_0^Tds  
\frac{1}{N^{d-1}} \sum_{\check x \in\T_N^{d-1}}
G_s(1, \frac {\check x} N ) \Big\{  \frac{1}{2cN}
\sum_{x_1=N(1-a-2c)+1}^{N(1-a)} \Big( \eta_s^{aN} (x_1,\check x)
-\eta_s^{\ell}(N-\ell,\check x)\Big)
\Big\}\Big | \Big] 
\; + \; R(N,a,c, G)\; ,
$$
where for $x_1\in[-N,N]$, $\check x =(x_2,\cdots,x_d)\in\T_N^{d-1}$ the
vector $(x_1,\check x)$ stands for the element $(x_1,x_2,\cdots,x_d)\in\L_N$. 
We denoted   by $ R(N,a,c, G)$  a quantity  
so that for  $G \in C^{1,2} \left ([0,T] \times
\L \right ) $, 
 $$
 \limsup_{c \to 0} \limsup_{a \to 0}\limsup_{N \to \infty}   |R(N,a,c,G)| =0.\;  
 \Eq(reste) 
 $$
The next step consists in replacing the density average over a small
macroscopic box of length $aN$ by a large microscopic box. More
precisely, for $N$ large enough the expectation of the last quantity
is bounded above by
$$
C \| G\|_\infty  \sup_{2\ell<|y|\le 2Nc}  \Es_{\mu^N}\Big[ \int_0^Tds  
\frac{1}{N^{d-1}} \sum_{\check x \in\T_N^{d-1}}
\Big|   \eta_s^{\ell}\big( (N-\ell,\check x)+y\big)
-\eta_s^{(\ell)}(N-\ell,\check x) \Big| \Big]
\; + \; R(N,a,c,\ell)\; ,
\Eq(bord5)$$
where for all $\ell$, $R(N,a,c,\ell)$ satisfy \equ(reste) and $C$ is a
positive constant.
Observe that the first term of the previous formula is not depending on $a$ but
only on $c,N$ and $\ell$.

 In view of the estimate \equ(bound-l1) and
Lemma \equ(dirichlet-l2) on the Dirichlet form $\DD_N^{0}$ and the
entropy, by the usual two blocks estimate, the first term of \equ(bord5) converges
to 0 an $N\uparrow\infty$, $c\downarrow 0$ and $\ell
\uparrow\infty$. That concludes the proof of Lemma \equ(bord-l2).
\qed


\medskip
\noindent {\bf Proof of Lemma \equ(bord-l3). }
The summation in \equ(bord4) contains two similar terms, we consider
the one corresponding to the summation of the right hand side of
$\L_N$. It is easy to see that the expectation in \equ(bord4) is
bounded above by
$$
\| G\|_\infty \frac{1}{N^{d-1}}\sum_{y\in\G_N^+} 
\Es_{\mu^N}\Big[ \int_0^Tds \Big| 
\eta_s^\ell(y-\ell e_1)-b(y/N) \Big|\Big]\; .
\Eq(bright)
$$
For any  fixed positive integer $\ell$  denote  by $\Gamma_0^\ell=\{(0,\hat x)\ : \  \ \hat x\in\T_N^{d-1},\ \ |\hat x|\le \ell \}
=(\{0\}\times \T_N^{d-1})\cap \Lambda_\ell(0)$, for notation see   \equ(cube.1).  
For  $u\in \Gamma$,  denote   $$
{\tilde D}_{\ell,0}^{b,u}
\big({f},  \nu  \big) 
\; =\; \frac 12 \sum_{\scriptstyle x\in
  \Gamma_0^\ell } 
\int {\tilde C}_0^{b} (u,x,\eta) \big( f  (\eta^x) -f  (\eta)\big)^2
d \nu (\eta) \; ,
$$
where  
$$
{\tilde C}_0^{b} (u,x,\eta)=
\h (x) \exp\Big\{ -\frac{\a (x) +\l_0(b (u)) }{2}\Big\} +(1-\h (x)) 
 \exp\Big\{ \frac{\a(x) +\l_0(b(u) )}{2}\Big\} \, . \Eq (rate2)
$$
The difference with the rate in \eqv (CD1) is  that here  $u$ is fixed. 
Let  $\nu_{b(u)}^{\a, N}$ be  the product   measure, see \eqv (inv1), where $\rho (\frac x N) \equiv b(u)$ for $\forall x\in \L_N$ and $\nu_{b(u)}^{\a, \ell}$ the  restriction of $\nu_{b(u)}^{\a, N}$    to $\{0,1\}^{\Lambda_\ell(0)}$.  
Let  $f: \SS_N \to \R$,  
  denote by $f^\ell$ the conditional
expectation of $f$  with respect to the $\sigma$-algebra generated by $\{\eta(z) \ : \ z\in \Lambda_\ell(0)\}$ :
$$
f^\ell(\xi)=\frac{1}{\nu_{ b(u)}^{\a, \ell}(\xi)}
\int \1_{\{\eta; \ \eta(z)=\xi(z),\ z\in \Lambda_\ell(0)\}}
f(\eta)d\nu_{b(u)}^{\a, N}(\eta)
\quad \text{for\ all}\quad \xi\in \{0,1\}^{\Lambda_\ell(0)}.
$$
Note that  
$\Big|\eta^\ell(0)-b(u)\Big|$ depends only on coordinates on the box $\Lambda_\ell(0)$, then by Fubini's Theorem, 
$$
\Es_{\mu^N}\Big[ \int_0^Tds \Big| 
\eta_s^\ell(y-\ell e_1)-b(y/N) \Big|\Big]=T\int \Big| 
\eta^\ell(0)-b(y/N) \Big|\big(\tau_{-(y-\ell e_1)}{\bar f}_T^{y,N}\big)^\ell (\eta)d\nu_{b({y\over N})}^{\a, \ell}(\eta)\; \Eq (M.1)
$$
where ${\bar f}_{T}^{y,N}=\frac1T \int_0^T f_s^{y,N} ds$
and for all $0\le s\le T$, $f_s^{y,N}$ is the density of $\mu_s^N$ with respect to
the product measure $\nu_{b({y\over N})}^{\a, N}$ with constant profile
$b({y\over N})$. The density $\big(\tau_{-(y-\ell e_1)}{\bar f}_T^{y,N}\big)^\ell$
stands for the conditional
expectation of $\tau_{-(y-\ell e_1)}{\bar f}_T^{y,N}$  with respect to the $\sigma$-algebra generated by $\{\eta(z) \ : \ z\in \Lambda_\ell(0)\}$.

Remark that, since the Dirichlet form is convex and since the conditional expectation is an average,
$$
\eqalign{
{\tilde D}_{\ell,0}^{b,\frac{y}{N}}
\Big(\sqrt { \big(\tau_{-(y-\ell e_1)}{\bar f}_T^{y,N}\big)^\ell}, \nu_{b(y/N)}^{\a, \ell}   \Big) 
&\le {\tilde D}_{\ell,0}^{b,\frac{y}{N}}
\Big(\sqrt { \tau_{-(y-\ell e_1)}{\bar f}_T^{y,N}}, \nu_{b(y/ N)}^{\a, N} \Big) 
\cr
\ &=\DD_{\ell,y-\ell e_1}^{b}
\Big(\sqrt { {\bar f}_T^{y,N}}, \nu_{b(y/N)}^{\a, N} \Big) \cr
\ &\le
\frac1T \int_0^T  \DD_{\ell,y-\ell e_1}^{b}
\Big(\sqrt { f_s^{y,N}}, \nu_{b(y/N)}^{\a, N} \Big)ds \; . 
}
\Eq(ddb)
$$
Applying  
  Lemma  \equ(dirichlet-l2b) 
we obtain from \equ(ddb)
$$\eqalign{
&N^{1-d}\sum_{y\in\Gamma_N}
{\tilde D}_{\ell,0}^{b,\frac{y}{N}}
\Big(\sqrt { \big(\tau_{-(y-\ell e_1)}{\bar f}_T^{y,N}\big)^\ell}, \nu_{b(y/N)}^{\a, \ell}   \Big) \cr
& \qquad\qquad\qquad
\le \frac1T \int_0^T \Big\{
N^{1-d}\sum_{y\in\Gamma_N}
\DD_{\ell,y-\ell e_1}^{b}
\Big(\sqrt { f_s^{y,N}}, \nu_{b(y/N)}^{\a, N} \Big)\Big\}ds\cr
&\qquad\qquad\qquad
\le 2\frac1T \int_0^T \Big\{
N^{1-d}\sum_{y\in\Gamma_N}
\DD_{\ell,y-\ell e_1}^{b}
\Big(\sqrt { h_s^{N}}, \nu_{\gamma(\cdot)}^{\a, N} \Big)\Big\}ds
+C_0\frac{\ell^{d+1}}{N^2}\cr
&\qquad\qquad\qquad
\le \frac{C_T}{N} +C_0\frac{\ell^{d+1}}{N^2}\; ,  
}\Eq (dd8)
$$
for some constant $C_T$ that depends on $T$.
By the same argument we obtain the bound on the Dirichlet form $\DD^0_{\ell,0}$,
$$
N^{1-d}\sum_{y\in\Gamma_N}
\DD^0_{\ell,0}\Big(\sqrt { \big(\tau_{-(y-\ell e_1)}{\bar f}_T^{y,N}\big)^\ell}, \nu_{b(y/N)}^{\a, \ell}   \Big) \le
\frac{C_T}{N} +C_0\frac{\ell^{d}}{N^2}\; .
\Eq(dd0)$$
 For $N$ fixed and large enough, there exists a constant $C_T$, such that for all positive integer $k\ge 1$,  applying \eqv (dd8) and \eqv  (dd0), we can bound by above  the expectation \equ(bright) as following 
$$\eqalign{
&T\| G\|_\infty N^{1-d}\sum_{y\in\G_{N}^+} 
\Big\{  \int \Big| 
\eta^\ell(0)-b(y/N) \Big| \big(\tau_{-(y-\ell e_1)}{\bar f}_T^{y,N}\big)^\ell
d\nu_{b(y/N)}^{\a, \ell}(\eta)
-k\; \DD^0_{\ell,0}\Big(\sqrt { \big(\tau_{-(y-\ell e_1)}{\bar f}_T^{y,N}\big)^\ell}, \nu_{b(y/N)}^{\a, \ell}   \Big) 
\cr
&\qquad\qquad\qquad\qquad\qquad\qquad\qquad\qquad\qquad
-\; k\; 
{\tilde D}_{\ell,0}^{b,\frac{y}{N}}
\Big(\sqrt { \big(\tau_{-(y-\ell e_1)}{\bar f}_T^{y,N}\big)^\ell}, \nu_{b(y/N)}^{\a, \ell}   \Big) \Big\}  \; +\frac{k}{N}\big(C_T +\frac{\ell^d(\ell+1)}{N}\big)\; .
}$$
This last expression is bounded above by

$$\eqalign{
&T\| G\|_\infty N^{1-d}\sum_{y\in\G_{N}^+} \sup_{f\in \AA_{\ell}^{\frac{y}{N}} }
\Big\{  \int \Big| 
\eta^\ell(0)-b(y/N) \Big| f(\eta)
d\nu_{b(y/N)}^{\a, \ell}(\eta)
-k\; \DD^0_{\ell,0}\Big(\sqrt {f}, \nu_{b(y/N)}^{\a, \ell}   \Big) 
\cr
&\qquad\qquad\qquad\qquad\qquad\qquad\qquad\qquad\qquad
-\; k\; 
{\tilde D}_{\ell,0}^{b,\frac{y}{N}}
\Big(\sqrt {f}, \nu_{b(y/N)}^{\a, \ell}   \Big) \Big\}  \; +\frac{k}{N}\big(C_T +\frac{\ell^d(\ell+1)}{N}\big)\; ,
}
\Eq(dd1)$$
where, for $u\in\Gamma$,
$$
\AA_{\ell}^u=\Big\{f:\ \ f\ge 0,\ \ \int f (\xi) d\nu_{b(u)}^{\a, \ell}(\xi)=1\Big\}\; .
$$
Further, since the function
$$
u\to \sup_{f\in \AA_{\ell}^{u} }
\Big\{  \int \Big| 
\eta^\ell(0)-b(u) \Big| f(\eta)
d\nu_{b(u)}^{\a, \ell}(\eta)
-k\; \DD^0_{\ell,0}\Big(\sqrt {f}, \nu_{b(u)}^{\a, \ell}   \Big) 
-\; k\; 
{\tilde D}_{\ell,0}^{b,u}
\Big(\sqrt {f}, \nu_{b(u)}^{\a, \ell}   \Big) \Big\}
$$
is continuous on $\Gamma$,  from Lemma \equ(lem-ergodic), for all positive integers $\ell$ and $k$, the limit when $N\uparrow \infty$ of the expression 
\equ(dd1) is equal to
$$
T\| G\|_\infty \int_\Gamma du \; \E \Big[
\sup_{f\in \AA_{\ell}^{u} }
\Big\{  \int \Big| 
\eta^\ell(0)-b(u) \Big| f
d\nu_{b(u)}^{\a, \ell}(\eta)
-k\; \DD^0_{\ell,0}\Big(\sqrt {f}, \nu_{b(u)}^{\a, \ell}   \Big) 
-\; k\; 
{\tilde D}_{\ell,0}^{b,u}
\Big(\sqrt {f}, \nu_{b(u)}^{\a, \ell}   \Big) \Big\}\Big]\; .
$$

Since $\int \Big|\eta_s^\ell(0)-b(u) \Big| f
d\nu_{b(u)}^{\a, \ell}(\eta)\le C_b$ for some positive constant $C_b$ that depends on $\|b\|_\infty$, the integral on $\Gamma$ in the last expression is bounded by
$$
\int_\Gamma du \; \E \Big[
\sup_{f\in \AA_{\ell, k, C_b}^u}
\Big\{ \int \Big| 
\eta^\ell(0)-b(u) \Big| f(\eta)d\nu_{b(u)}^{\a, \ell}(\eta)
\Big\}\Big]\; ,
$$
where for a positive constant $C$, $\AA_{\ell, k, C}^u$ is the following set of densities,
$$
\AA_{\ell, k, C}^u=\Big\{f\in \AA_{\ell}^u\; ,\quad
 {\tilde D}_{\ell,0}^{b,u}
\Big(\sqrt {f}, \nu_{b(u)}^{\a, \ell}   \Big) \le \frac{C}{k}\; ,
\quad
\DD^0_{\ell,0}\Big(\sqrt {f}, \nu_{b(u)}^{\a, \ell}   \Big) \le \frac{C}{k}
\Big\}\; .
$$
We first consider the limit when $k\uparrow\infty$ and use the usual technics in the replacement lemma.
Since for any $\ell >1$, any constant $C>0$ and any $u\in\Gamma$ 
the sets $\AA_{\ell, k, C}^u$ are compacts for the weak topology, for all $\ell >1$
$$
\limsup_{k\to \infty} 
\sup_{f\in \AA_{\ell, k, C}^u}
\left\{ \int \Big| 
\eta^\ell(0)-b(u) \Big| f(\eta)d\nu_{b(u)}^{\a, \ell}(\eta)
\right\}=
\sup_{f\in \AA_{\ell,  C}^u}
\left\{ \int \Big| 
\eta^\ell(0)-b(u) \Big| f(\eta)d\nu_{b(u)}^{\a, \ell}(\eta)
\right\}\; ,
$$
where
$$
\AA_{\ell,   C}^u=\Big\{f\in \AA_{\ell}^u\; ,\quad
 {\tilde D}_{\ell,0}^{b,u}
\Big(\sqrt {f}, \nu_{b(u)}^{\a, \ell}   \Big) =0\; ,
\quad
\DD^0_{\ell,0}\Big(\sqrt {f}, \nu_{b(u)}^{\a, \ell}   \Big)=0
\Big\}\; .
$$
By  dominated convergence theorem, it is then enough to show that,
$$
\limsup_{\ell\to \infty} \E \Big[
\sup_{f\in \AA_{\ell, C}^u}
\left\{ \int \Big| 
\eta^\ell(0)-b(u) \Big| f(\eta)d\nu_{b(u)}^{\a, \ell}(\eta)
\right\}
\Big]=0\; .
$$
Now, it is easy to see that, due to the presence of the jumps of particles in
the Dirichlet form $\DD^0_{\ell,0}$  and the presence of
 the creation and destruction of particles in ${\tilde D}_{\ell,0}^{b,u}$
 the set $\AA_{\ell, C}^u=\{1\}$. Thus, to conclude the proof of the lemma, it remains to apply the usual law of large numbers.
\qed


\medskip
\noindent{\bf Proof of Proposition \equ(lem2b)}
Let $Q^*$ be a limit point of the sequence $(Q_{\mu^N})_{N\ge 1}$ and
let $(Q_{\mu^{N_k}})_{k\ge 1}$ be a sub-sequence converging to
$Q^*$. 
By Lemma \equ(lemIB1) $Q^*$ is concentrated on the
trajectories that are in $L^2([0,T];H^1(\L))$.
For $0<c<1$ and for $\mu(\cdot,\cdot) \in D([0,T],\MM_1^0(\L))$, such that 
$\mu(t,du)=\rho(t,u)du$ with $\rho(\cdot,\cdot)\in L^2([0,T];H^1(\L))$, 
denote by $F_{c}^{G}(\mu )$ the functional
$$\eqalign{
 F_{c}^{G}\big(\mu (\cdot, \cdot) \big)&
=\int_0^T  ds   \int_{\L_{(1-c)}}  du \Big\{  G_s (u)\;   {(2c)}^{-1}
\Big[ \rho(s,u+{c}e_1) -  \rho(s,u-{c}e_1)\Big]  \Big\}\cr
\ &\quad
+\int_0^T  ds   \int_{\L} du \partial_{e_1}G_s(u)  \rho(s,u) 
-\int_0^T ds\Big\{  \int_\Gamma b(u) {\hbox{\bf
n}}_1(u)  G_s (u) 
 \hbox{d} \hbox{S} \Big\}\; .
}$$  
From Lemma \equ(bord) and the continuity of the function
$\mu\to \hat F_{a,c}^{G}(\mu )$, we have
$$
\limsup_{c\to 0}E^{Q^*}
\Big[ \Big| F_{c}^{G}(\mu) \Big| \Big] =0\; . \Eq (M1)
$$
On the other hand, an integration by parts and Taylor expansion up  to the
second order of the function $G_s (\cdot)$  permit to rewrite $F_c^G$ as
$$\eqalign{
 F_{c}^{G}\big(\mu (\cdot, \cdot) \big)&
=\int_0^T   \frac{1}{2c}  \int_{(\L\setminus\L_{(1-2c)})^+} G_s(u)
\rho(s,u) du ds
- \int_0^T   \frac{1}{2c}  \int_{(\L\setminus\L_{(1-2c)})^-} G_s(u)
\rho(s,u) du ds\cr
\ &\quad
-\int_0^T ds  \int_\Gamma b(u) {\hbox{\bf
n}}_1(u)  G_s (u) 
\hbox{d} \hbox{S}  \; +\; R(c)\; ,
}$$
where $R(c)\equiv R(G,c)$ is a function vanishing as $c\downarrow 0$.
 Further one has, see Theorem 5.3.2. of [EG],  that  
 $$ \lim_{r \to 0}  \frac 1 { |B (u, r) \cap \Lambda| }   \int_{ B (u, r) \cap \Lambda} \rho(s,y) dy 
 = \Tr(\rho(s,u))  \qquad \hbox {a.e} \quad   u \in \G, \forall s,    \Eq (M5)$$ 
and then  by dominated convergence theorem 
$$
\lim_{c\to 0}  F_{c}^{G} \big(\mu (\cdot, \cdot) \big)=
 \int_0^T ds  \int_\Gamma \Big(\Tr(\rho(s,u) )- b(u)\Big) {\hbox{\bf
n}}_1(u)  G_s (u) 
\hbox{d} \hbox{S}  \; .
$$
This together with  \eqv (M1) implies 
$$
E^{Q^*}\Big[  \Big|\int_0^T ds  \int_\Gamma \Big(\Tr(\rho(s,u) )- b(u)\Big) {\hbox{\bf
n}}_1(u)  G_s (u) 
\hbox{d} \hbox{S} \Big|\Big]=0\; ,
$$
which  concludes the proof.
\qed

\vskip0.5cm 
\chap{4  Proof of Theorem \eqv(hydro1) }4
\numsec= 4
\numfor= 1
\numtheo=1 
  
The main problem in proving  Theorem \equ(hydro1) is  that  we cannot associate  to the stationary measure $\mss$ a macroscopic profile according to definition \eqv (o.1).     If this would be the case  the result  would  be a corollary of   Theorem \eqv (th-hydro).
 Denote by $ \Qss^N:= Q^{N,\a}_{\mss}$ the probability measure on the Skorohod space
$D\big([0,T], {\cal M}\big)$ induced by the Markov process $(\pi_t^N)\equiv (\pi_N(\eta_t))$, when the initial measure is $\mss$.
Denote by ${\cal A}_T\subset D\big([0,T], {\cal M}\big)$ the class of
profiles $\rho(\cdot,\cdot)$ that satisfies conditions (IB1), (IB2) and (IB3). 
The first step  to show  Theorem \eqv  (hydro1) consists in proving that all limit points of the sequence
$(\Qss^N)$ are concentrated on  ${\cal A}_T$:

\medskip
\noindent{\bf  \Proposition (pr-hydrost)} 
{\it The sequence of probability measures
$(\Qss^N)$ is weakly relatively compact and  all its converging
subsequences converge to the some limit $\Qss^*$ that is concentrated on
the absolutely continuous measures $\pi(t,du)=\rho(t,u)du$ whose density
$\rho$ satisfying (IB1), (IB2) and (IB3).
} 
\medskip

The proof   of  Proposition \eqv  (pr-hydrost) follows the same steps  needed to   show Theorem \equ(th-hydro).
We just have to show the analogous of Lemmas \equ(dirichlet), \equ(dirichlet-l1), 
\equ(dirichlet-l2) and \equ(dirichlet-l2b) when the  measure  $\mu^N$ in the statements of these lemmas is replaced by   $\mss$.
 The only  lemma  to be  slightly modified  is  Lemma  \equ(dirichlet-l2), see  Lemma \eqv  (dirichlet-lss) given next.      
 Recall that $\gamma\colon\Lambda\to (0,1)$ is a smooth profile equal to
$b$ at the boundary of $\Lambda$. Let $h^N$ be the density
of $\mss$ with respect to the  measure $\nu_{\gamma(\cdot)}^{\a,N}$. 

\medskip
\noindent {\bf \Lemma (dirichlet-lss)} 
{\it There exists positive constant $C=C(\|\nabla \gamma\|_\infty)$ depending only  on $\gamma(\cdot
)$ such that  for any $a>0$
$$
  (1-a)  
\DD_N^{0}  (\sqrt {h^{N}}, \nu_{\gamma(\cdot)}^{\a,N} )
 + 
\DD_N^{b}  (\sqrt {h^{N}}, \nu_{\gamma(\cdot)}^{\a,N} )\; \le\;
\frac{C}{a} N^{d-2}\; .
$$
}

\medskip

\noindent{\bf Proof.} 
 By the stationary of $\mss$,
$$
 \partial_t H_N(t)= \int_{\SS_N} h^N \LL_N \log \big( h^N  \big) d \nu_{\gamma(\cdot)}^{\a,N} =0\; . 
$$
Recalling that the generator $\LL_N$ has two pieces and applying 
  the basic inequality
$
a\big(\log b -\log a  \big)\le -\big( \sqrt{a}-\sqrt{b}\big)^2+\big(b-a\big)
$
for positive $a$ and $b$, we obtain
$$\eqalign{
 0=\int_{\SS_N} h^N \LL_N \log \big( h^N  \big) d \nu_{\gamma(\cdot)}^{\a,N} & \; \le\;
-2N^2 \DD_N^{0} \big(\sqrt {h^{N}},\nu_{\gamma(\cdot)}^{\a,N}  \big) 
-2N^2 \DD_N^{b} \big(\sqrt {h^{N}},\nu_{\gamma(\cdot)}^{\a,N}   \big) \cr
\ &\ \ \ 
+N^2\int_{\SS_N}  \LL^0_{N} h^N  d \nu_{\gamma(\cdot)}^{\a,N}  
+ N^2\int_{\SS_N} \LL^b_{N} h^N  d \nu_{\gamma(\cdot)}^{\a,N}  \; .
}
$$
We then apply the  same computation as in the proof of Lemma \equ(dirichlet-l2),
(\equ(entr1) and \equ(entr2)).
\qed
 \medskip
\noindent{\bf Proof of Theorem  \equ(hydro1)}
\medskip 
Let $\Qss^*$
be a limit point of $ (\Qss^N) $ and   $(\Qss^{N_k})$ be a
sub-sequence converging to $\Qss^*$.   Let $\bar \r$  be the stationary solution  of \eqv (heq1), see \eqv (SS.1).  We have  
by   Proposition \eqv  (pr-hydrost) the following: 
$$\eqalign{
\lim_{k\to \infty}q_{N_k} (\a) & \;: =\; \lim_{k\to \infty}
\Qss^{N_k} \Big\{ \Big|\big< \pi_T^N,G \big>- \big< \bro(u) du,G\big>
\Big|\Big\}\cr
\  & \; =\; 
\Qss^* \Big\{ \Big|\big< \rho(T,\cdot) ,G\big>- \big< \bro(u) du,G\big>
\Big| \1 \{{\cal A}_T\} \big(\rho\big)   \Big\}\cr
 &\ \ \; \le \; \| G\|_\infty
\Qss^* \Big\{ \big\| \rho(T,\cdot) -\bro(\cdot) \big\|_1
\1 \{{\cal A}_T\} \big(\rho\big)\}   \Big\}\; .
}
$$
Denote by  $\rho^0(\cdot,\cdot)$ (resp.
$\rho^1(\cdot,\cdot)$)  the element of ${\cal A}_T$ with initial condition $\rho^0(0,\cdot)\equiv 0$ 
(resp. $\rho^1(0,\cdot)\equiv 1$). From Lemma \equ(S5), each profile
$\rho(\cdot,\cdot)\in {\cal A}_T$ is such that for all $t\ge 0$,
$\lambda\Big\{u\in\L\ : 0 \le  \rho^0(t,u) \le     \r(t,u) \le \rho^1(t,u) \le 1\Big\}=1$ and
$\lambda\Big\{u\in\L\ :  \rho^0(t,u) \le    \bro(u) \le \rho^1(t,u)\Big\}=1$,
where $\lambda$ is the Lebesgue measure on $\L$. Therefore
$$
\lim_{k\to \infty} q_{N_k}(\a) \;\le\; \| G\|_\infty \; 
\big\| \rho^0(T,\cdot) -\rho^1(T,\cdot) \big\|_1\;,  \quad \P=1.
$$
To conclude the proof,  it is enough to let $T\uparrow \infty$ and to apply Theorem
\equ(AA).
\qed

\newpage 
  \chap {5   Appendix  } 5
\numsec=5
\numfor=1
\numtheo=1
In this section we show the global stability of the stationary solution of \eqv (heq1).   
 
\medskip \noindent 
  {\bf \Theorem(AA)} {\bf  Global stability}. {\it  Let $D (\cdot)$ be  Lipschitz. 
  Let $\r(t, \r_0)$ be the  solution  of \eqv (heq1) with initial datum $\r_0$,  $0\le \r_0 (u) \le 1$, $u \in \L$,  and   
   $\bar \r$   the stationary solution of \eqv (heq1). We have 
  $$ \lim_{t \to \infty}\int_{\L}  |\r(t, u)- \bar \r (u)|^p du =0 $$
  for all $p \ge 1$ . } 
  \vskip0.5cm \noindent 
   The proof   of the theorem  is based on an extensive use of    monotone methods, see [S].  We were not able to find the precise reference, so  we   briefly   sketch it  for completeness. 
  We need to introduce some extra notation.  Let  
$ {\CC}^{1,2}\big([0,T]\times\L \big)$ be    
the space  of functions  from 
  $[0,T]\times \L$ to $\R$  twice continuously differentiable in $\L$
with continuous time derivative.  Denote by 
 $$ \GG:= \left \{   G \in  {\CC}^{1,2}\big([0,T]\times\L \big),  G(t,u)=G_t(u)\quad \hbox { pointwise  positive},  \quad  G(t, u )=0, \forall u \in \G, \forall  t \in [0,T] \right \}. $$
  It is convenient to reformulate  the notion of  weak solution of \eqv (heq1)  as following. 
A   function $\rho(\cdot, \cdot ):[0,T]\times \L\to [0,1] $    is a weak solution   
of the initial-boundary value problem  \eqv (heq1) if    $ \r \in  L^2\big(0,T;H^1(\Lambda)\big)$ and   for every $ G \in \GG$
$$\eqalign{
& \int_\L du \big\{ G_T(u)\rho(T,u)-G_0(u)\rho_0(u)\big\} -
\int_0^T ds \int_\L d u \,  (\partial_s G_s)(u)\rho(s,u) \cr
& \quad =\;   \sum_{i, j} \int_0^T d s \Big\{\int_\L d u \, A_{i,j}(\rho(s,u))\frac { \partial^2} {\partial_{i, j} }   G_s(u)  - \int_\G  A_{i,j}(b(u)) \partial_{n_1} G(s,u) d S    \Big \}
} \Eq (B1)$$
where 
$  A_{i,j} (\rho) = \int_0^{\r} D_{i,j} (\r' )  d \r' $.
 A   function $\rho^{+}(\cdot, \cdot ):[0,T]\times \L\to \R $ is  
a weak  {\bf  upper solution  } 
of the initial-boundary value problem  \eqv (heq1)
if  $ \r^+ \in  L^2\big(0,T;H^1(\Lambda)\big)$  and for  all $ G \in \GG$
    we have 
$$  \left \{  \eqalign{
& \quad   \;   \sum_{i, j} \int_0^T d s \Big\{\int_\L d u \, A_{i,j}(\rho^+(s,u))\frac { \partial^2} {\partial_{i, j} }   G_s(u)  - \int_\G  A_{i,j}(\rho^+(s,u)) \partial_{n_1} G(s,u) d S       \Big \}  \cr & - \int_\L du \big\{ G_T(u)\rho^+(T,u)-G_0(u)\rho^+_0(u)\big\} -
\int_0^T ds \int_\L d u \,  (\partial_s G_s)(u)\rho^+(s,u) 
  \le 0,  \cr & 
\Tr(\rho^{+} ( t, \cdot)) \ge b(\cdot)  \quad \hbox {on} \quad  \G  \cr &
 \rho^{+} ( 0, u) \ge \r_0(u) \quad u \in \L 
}  \right.  \Eq (B3)$$
 A    weak  {\bf  lower solution  }  $\rho^{-}(\cdot, \cdot ):[0,T]\times \L\to \R $ is  
 defined  reversing the inequality in \eqv  (B3).

 \noindent
 By a solution of the stationary problem   \eqv  (heq1)  we mean a function $  \bar \r  \in H^1 (\L)$ so that  for  all $G \in  {\CC}^{2}\big(  \L \big)$,  {\bf {pointwise  positive}}  vanishing on $\G$
$$ \sum_{i, j}  \Big\{\int_\L d u \, A_{i,j}(\bar \r (u)))\frac { \partial^2} {\partial_{i, j} }   G (u)  -  \int_\G  A_{i,j}(b(u)) \partial_{n_1} G(u) d S \Big \}   =0 
\Eq (B2)$$
As before we define   upper  and  lower  solutions of the stationary problem \eqv (B2).     
 A function  $ {\bar \r}^+ $ is an  upper solution for the stationary problem \eqv (B2) if 
 $ {\bar \r}^+   \in H^1 (\L)$ and 
  for  all $G \in  {\CC}^{2}\big(  \L \big)$,  pointwise  positive  vanishing on $\G $, 
$$ \left \{   \eqalign { &
 \sum_{i, j}  \Big\{\int_\L d u \, A_{i,j}( {\bar \r}^+(u)))\frac { \partial^2} {\partial_{i, j} }   G (u)  -  \int_\G  A_{i,j}( {\bar \r}^+ (u)) \partial_{n_1} G(u) d S \Big \}    \le 0 \; , \cr & 
 \Tr( {\bar \r}^+) \ge b \quad \hbox {on}\quad \G ,} \right.   \Eq (SS8)
$$
  A lower  solution of the stationary problem \eqv (B2) is defined  reversing the inequality in \eqv (SS8).
\medskip

To  apply the monotone  method we   first  show the following comparison principle.
\medskip
\noindent {\bf \Lemma (stable-region)} 
{\it Let $\rho^1$ (resp. $\rho^2$) be   a lower solution (resp.  upper solution) of \equ(heq1),  
$\partial_t \rho^i  \in L^2\big(0,T;H^{-1}(\Lambda)\big)$, for $i=1,2$.
If there exists $s\ge 0$ such that
$$\lambda \big\{ u\in  \L\ :\ \ 
\rho^1(s,u) \le\rho^2(s,u)  \big\}=1\; , $$ 
where $\lambda$ is the Lebesgue measure on $ \L$,
then for all $t\ge s$
$$\lambda \big\{ u\in  \L\ :\ \ 
\rho^1(t,u) \le\rho^2(t,u)  \big\}=1.
$$
}\vskip0.5cm 
\noindent{\bf Proof}
Take    $s< t<T$ and 
  $\delta >0$.  Denote by $F_\delta$ the function defined by
$$ 
F_\delta (a):=\; \frac{a^2}{2\delta }\1 _{\{0\le a\le \delta \}}
\; +\; \big(  a -\delta/2\big)
 \1_{ \{a >\delta\}}, \quad    a \in \R. 
$$
Let $A_\delta:=A_\delta(T)$ be the set
$$A_\delta  =\Big\{ (t,u)\in
[0,T]\times  \L\ :\ 0\le \rho^1(t,u) - \rho^2(t,u)\le
\delta \Big\}.
$$ 
By  definition    $\Tr(\rho^1-\rho^2)\le 0 $ a.e. and  therefore $\Tr\big(F_\delta'(\rho^1-\rho^2)\big)=0$.
 Since $\rho^1 $  ( $\rho^2 $ )is   lower  (upper) solution    of \equ(heq1),   we  have that
 $$
\eqalign{ \int_s^t d\t \frac {\partial } {\partial \tau} \int_\L  F_\delta \Big (\rho^1(\t,u) -\rho^2(\t,u)\Big)=
&\int_ \L du \, F_\delta \Big (\rho^1(t,u) -\rho^2(t,u) \Big) 
-\int_ \L du \, F_\delta \Big (\rho^1(s,u) -\rho^2(s,u) \Big) \cr
&\quad \le \; - \delta^{-1} \int_s^t d\tau \int_{A_\delta} du\,
\nabla (\rho^1 -\rho^2) \cdot \Big\{ D(\rho^1)  \nabla
\rho^1 -D(\rho^2)  \nabla \rho^2\Big\} \cr 
&\qquad = \; - \delta^{-1} \int_s^t d\tau  \int_{A_\delta} du\,
\nabla (\rho^1 -\rho^2)\cdot D(\rho^1)  \nabla (\rho^1 -\rho^2) \cr
&\qquad\quad -\; \delta^{-1}\int_s^t d\tau \int_{A_\delta} du\,
\nabla (\rho^1 -\rho^2) \cdot \big\{ D(\rho^1)- D(\rho^2) \big\}
\nabla \rho^2 \; . }  \Eq (T1)
$$ 
Since   $D(\cdot)$ is strictly positive,  see \eqv  (elliptic), 
the third line  of \eqv (T1) can be estimated  by  above 
$$
  -\frac 1 \d \int_s^t d\tau  \int_{A_\delta} du\,
\nabla (\rho^1 -\rho^2)\cdot D(\rho^1)  \nabla (\rho^1 -\rho^2)    \le  - \frac 1 {\d C }  \int_s^t d\tau  \int_{A_\delta} du\, \Vert \nabla
(\rho^1 -\rho^2) \Vert^2.  \Eq (T2)  
$$
 Further, by the  Lipschitz property of $D(\cdot)$  we have on the set $A_\delta$,
$\sup_{1\le i,j\le d}| D_{i,j}(\rho^1)- D_{i,j}(\rho^2)| \le M| \rho^1 -\rho^2| \le M\delta$ for some
positive constant $M$. 
By Schartz inequality, the last line of \eqv (T1)   is bounded by
$$ 
\delta^{-1} M A \int_s^t d\tau  \int_{A_\delta} du\, \Vert \nabla
(\rho^1 -\rho^2) \Vert^2\; +\; 
\delta M A^{-1} \int_s^t d\tau  \int_{A_\delta} du\, 
\Vert \nabla \rho^2 \Vert^2 \Eq (T3)
$$ 
for every $A>0$.
By \eqv (T1), \eqv (T2),  \eqv (T3) and  choosing  $A=M^{-1} C^{-1}$  to cancel  the  term  in \eqv (T2) and the first term of  \eqv  (T3)  we have 
$$
\eqalign{ 
&\int_ \L du \, F_\delta \Big(\rho^1(t,u) -\rho^2(t,u) \Big)
-\int_ \L du \, F_\delta \Big (\rho^1(s,u) -\rho^2(s,u) \Big)\cr
&\qquad \le \; \delta C^{-1} M^2  \int_0^T d\tau  \int du\, 
\Vert \nabla \rho^2 \Vert^2 \; .
}
$$ 
  Letting $\delta\downarrow 0$, we conclude the proof of the lemma
because $F_\delta (\cdot)$ converges   to the function $F(a)=a\1_{a\ge 0}$ as
$\delta\downarrow 0$. 
\qed 

\medskip
By  Lemma \equ (stable-region)  we immediately obtain the following corollaries.

\noindent {\bf \Corollary (coro1)} {\it  Let $m_0 :\Lambda \to [0,1]$ be a measurable function. There is a unique weak solution  $\rho (t, m_0)$ of the equation \equ (heq1) with initial datum $m_0$. }


 \medskip \noindent {\bf \Corollary (S1)} 
{\it Let $ m_0$ be a lower stationary solution of  \eqv (B2).  Let $\rho (t, m_0)$ be the solution of \eqv (B1) with initial datum $m_0$  then $  \rho (t, u) \ge    m_0(u) $ a.e in $(u,t)$. } 
\medskip \noindent 
The proof is an immediate consequence of  Lemma   \equ (stable-region) with $\r^1:=m_0$ and $\r^2:=\r$.     
When the initial datum of solution of \eq (B1) is an upper  stationary solution we have:
  \medskip
\noindent {\bf \Corollary (S3)} 
{\it Let $ m_1$ be a  upper  stationary solution of  \eqv (B2). Let $\rho (t, m_1)$ be the solution of \eqv (B1) with initial datum $m_1$  then $  \rho (t, u) \le    m_1(u) $  for  $ t \in [0,T]$ and $u \in \L$. } 

Next we  show that  when    a lower (upper) stationary solution $ m_0$ ($m_1$) is  taken as initial datum,  
the corresponding solution $\r(t,m_0)$ ($\r(t,m_1)$) is monotone nondecreasing  (nonincreasing)  in time. 
\medskip \noindent {\bf \Lemma (S2)} 
{\it  Under the assumptions of Corollary \eqv (S1)  $\rho (t, m_0)$ is a nondecreasing  solution of \eqv (heq1) for all $ t \in [0,T]$. } 
\vskip0.3cm \noindent 
\proof    Corollary \eqv (S1)  implies that  $\r (s, m_0) \ge m_0$ for all $s\ge 0$, since $m_0$  lower solution.  Let     $ \r (t; \r (s, m_0)) $ be the solution of  \eqv (B1) starting at time $t=0$ from $ \r (s, m_0)$. 
Then   $ \r (t; \r (s, m_0))  \ge \r (t, m_0)$  since the initial datum  $\r (s, m_0)\ge m_0$.
But  $ \r (t; \r (s, m_0)) = \r (t+s, m_0)$ by uniqueness of weak solution  then  $\r (t+s, m_0) \ge  \r (t, m_0) \ge m_0$.  \qed 

 \medskip \noindent {\bf \Lemma (S4)} 
{\it  Under the assumptions of Corollary \eqv (S3)  $\rho (t, m_1)$ is a nonincreasing  solution of \eqv (heq1) for   $ t \in [0,T]$. } 
 \vskip0.3cm \noindent    The proof is similar to the one of Lemma \eqv (S2).


\noindent {\bf \Lemma (S5)} {\it Let  $m_0$  be a lower solution  and    $m_1  $ be an  upper solution of \eqv (B2), $ m_0(\cdot) \le m_1(\cdot)$ a.e in $\L$,  we have }
$$ m_0 \le \r (t; m_0) \le \r (t; m_1) \le m_1 \quad \forall t \in (0, \infty)$$
\vskip0.5cm \noindent 
The proof    is  an immediate consequence of  the previous results. \qed
 
 \noindent {\bf \Lemma (S6)} {\it  Under the assumption of Lemma \eqv (S5) 
 the solutions $ \r (t; m_0)$ and  $\r (t; m_1)$ exist  for all $ t \in [0, \infty)$ and
 they converge in   $ L^p  (\Lambda)$  for   $p \in [1, \infty)$  to limits $\rho_{\star}(\cdot)$ and $\rho^{\star} (\cdot)$, both solutions of  \eqv (B2).   Further 
 $$ \rho_{\star}(u) \le  \rho^{\star}(u)  \quad a.e.$$ }
 \vskip0.5cm \noindent 
 \proof Since  $\r (t; m_0)$ is nondecreasing in $t$ and $ \r (t; m_0) \le m_1$  for any $t\ge 0$,   $\r (t; m_0)$  converges almost everywhere in $\Lambda$ as $ t\to \infty$ and 
  $\rho_{\star} (\cdot) \in L^\infty ( \L)$.  By the monotone convergence theorem
  $\r (t; m_0) \to  \rho_{\star} (\cdot) $  for  $p \in [1, \infty)$. 
  Next we show that $ \rho_{\star} (\cdot)$ solves \eqv (B2). 
  Take  as test function in \eqv (B1) the following function
 $$\b(t) F(u); \quad F(u)>0; \quad  C \ge  \b(t)>\delta>0;  \quad \b'(t)\ge 0,  (u, t) \in \L\times  \R^+ $$
 $\b \in C^2(R^+)$,  $F \in C^2( \L)$ vanishing at the boundary.  
 Then  for all $t >0$, see  \eqv (B1), we have 
  $$\eqalign{  
& \int_\L du \big\{ \b(t) F(u)\rho(t,u)-\b(0) F(u)\rho_0(u)\big\} -
\int_0^t ds  \b'(s) \int_\L d u  F(u)\rho(s,u) \cr
& \quad =\;   \sum_{i, j} \int_0^t d s \b(s)  \Big\{\int_\L d u \, A_{i,j}(\rho(s,u))\frac { \partial^2} {\partial_{i, j} }  F(u)  - \int_\G  A_{i,j}(b(u)) \partial_{n_1} F(u) d S    \Big \}.
} \Eq (B3a)$$
Divide by $t $ the  left  and  right   side of \eqv (B3a)  and then let $t \to \infty$.
For the left side we have   
$$  \frac 1 t \left \{ \int_\L du \big\{ \b(t) F(u)\rho(t,u)-\b(0) F(u)\rho_0(u)\big\} -
\int_0^t ds  \b'(s) \int_\L d u  F(u)\rho(s,u) \right \}  \to 0. \Eq (B3b) $$
By continuity of $ A (\cdot)$ and  since by assumption  $\lim_{s \to \infty} \b (s) =\b(\infty)>0 $
$$ \eqalign { &  \lim_ {t \to \infty} \frac 1 t \sum_{i, j} \int_0^t d s \b(s)  \Big\{\int_\L d u \, A_{i,j}(\rho(s,u))\frac { \partial^2} {\partial_{i, j} }  F(u)-  \int_\G  A_{i,j}(b(u)) \partial_{n_1} F(u) d S  \Big \} \cr &
= \b(\infty) \sum_{i, j} \Big\{\int_\L d u \, A_{i,j}(\rho_\star (u) )\frac { \partial^2} {\partial_{i, j} }  F(u) -  \int_\G  A_{i,j}(b(u)) \partial_{n_1} F(u) d S  \Big \}. }  
    \Eq (B4)$$
  By \eqv   (B3b) we then obtain  $$ \b(\infty) \sum_{i, j} \Big\{\int_\L d u \, A_{i,j}(\rho_\star(u))\frac { \partial^2} {\partial_{i, j} }  F(u) -  \int_\G  A_{i,j}(b(u)) \partial_{n_1} F(u) d S  \Big \}=0.$$
  Therefore $\rho_\star$  is a solution of  \eqv (B2).
 The same can be argued for $\r^*$. 
   \qed  
  \vskip1.cm \noindent  
  The  proof of Theorem \eqv (AA) is a simple consequence of Lemma \eqv (S6) and   the unicity of the stationary solution 
$ \r^*=\r_{\star}$ of \eqv (heq1).   

 \medskip

 \vskip 0.5cm
 \noindent{\bf Acknowledgements}
Enza Orlandi  is indebted to Assunta Pozio  (La Sapienza, Roma) for helpful  suggestions  about  Section 5.  She further   
    thanks the University of Rouen and the Institut Henri Poincar\'e - Centre Emile Borel, (workshop 
     M\'ecanique statistique, probabilit\'es et syst\`emes de particules 2008)
  where part of the work has been done. 
Mustapha Mourragui  thanks  J. S. Farfan Vargas    for helpful discussions concerning the hydrostatic.   We thank   Claudio Landim for suggesting the approach
to prove hydrostatic

 \goodbreak
\vskip1truecm
\noindent {\bf References}
\item{[BSGJL]}    L. Bertini, A. De Sole, D. Gabrielli, G. Jona-Lasinio, C. Landim,   {\sl Large deviation approach to non equilibrium processes in stochastic lattice gases. }Bull. Braz.  Math. Soc. (N.S.) {\bf 37}   611-643  (2006). 

\item{[BSGJL]}   L. Bertini, A. De Sole, D. Gabrielli, G. Jona-Lasinio, C. Landim,   {\sl Large deviations of the empirical current in interacting particle systems.}   Teor. Veroyatn. Primen.  51  (2006),  no. 1, 144-170;  translation in  Theory Probab. Appl.  51  (2007),  no. 1, 2-27 

 \item{[DL]}   R. Dautray and J.L. Lions,   {\sl  Analyse math\'ematique et calcul num\'erique   pour le sciences et les techniques}.  tome 3, Masson, Paris    (1985)

\item{[DFIP]}    A. De Masi, P.  Ferrari, N. Ianiro, E.Presutti,   {\sl Small deviations from local equilibrium for a process which exhibits hydrodynamical behaviour.  II}. {  J. Stat. Phys. } {\bf 29},  81-93 (1982)

\item{[DEL]} B.  Derrida, C. Enaud, J.L. Lebowitz, 
{\sl The asymmetric exclusion process and Brownian excursions ,}  
{   J. Stat. Phys. } {\bf 115},  365-383 (2004)

\item{[DELO]} B.  Derrida, C. Enaud, C. Landim, S. Olla,
{\sl  Fluctuaction  in the weakly   asymmetric exclusion process with open boundary conditions.}  
{  J. Stat. Phys. } {\bf 121},  271-289, (2005)


  \item{[E]} L.C. Evans,
  {\sl  Partial Differential Equation,}  
  {American Mathematical Society}, (1998).  
 
 \item{[EG]} L.C. Evans and R. F. Gariepy,
 {\sl Measure Theory and Fine Properties of Functions,} 
 {Studies in Advanced Mathematics},  (1992). 

\item{[ELS1]}  G. Eyink, J. L. Lebowitz and H. Spohn,
Hydrodynamics of
 {\sl  Stationary Nonequilibrium States for Some Lattice Gas Models}. 
 Commun. Math. Phys.  {\bf 132},  252-283  (1990).

\item{[ELS2]}  G. Eyink, J. L. Lebowitz and H. Spohn,
  {\sl Lattice Gas Models in contact with Stochastic Reservoirs: Local Equilibrium and Relaxation to the
Steady State}.
  Commun. Math. Phys.  {\bf 140},   119--131 (1991).

 \item{[FM]} A. Faggionato  and F. Martinelli,
 {\sl Hydrodynamic limit of a disordered lattice gas}.
 {  Prob. Th. Rel. Fields} {\bf 127}, (4),  (2003), 535--608.

 \item{[FLM]}   J. S. Farfan Vargas, C. Landim and M. Mourragui,
 {\sl  Hydrodynamic behavior and large deviations of
  boundary driven exclusion processes in dimension $d>1$.}
in preparation

   \item{[GPV]} M.Z. Guo, G. Papanicolau   and   S.R.S.Varadhan,
 {\sl Nonlinear diffusion limit for a system with nearest neighbor interactions}. 
    Comm. Math. phys.  {\bf 118},   31--59, (1988).

\item{[KL]} C. Kipnis    and C. Landim, 
{\sl Hydrodynamic limit of interacting particle systems}. 
Springer-Verlag,  (1999).

\item{[LSU]} 
O. A. Ladyzenskaja, V.A. Solonnikov, N.N. Uralceva,
 {\sl Linear and quasi linear equations of parabolic type}. 
{  AMS}. {\bf 23},  (1968).
1998. 

\item{[KLO]} 
 C. Kipnis, C. Landim and S.Olla,
 {\sl Macroscopic properties of a stationary nonequilibrium distribution for a non-gradient interacting
particle system}. 
{ Ann. Inst. H. Poincar\'e}. {\bf 31}, 191--221, (1995). 

\item{[KW]} K.W. Kehr,  T. Wichman, 
 {\sl Diffusion Processes: experiment, theory of simulations }.
{Lectures Notes In Physics  } {\bf 438}, (1994)

\item{[Li]}T.  Liggett, 
 {\sl Interacting particles systems.} (1985),
 Springer, Berlin.

\item{[LMS]} C. Landim, M. Mourragui and S. Sellami,
  {\sl Hydrodynamical limit for a nongradient interacting particle  
system with stochastic reservoirs}.
{  Probab. Theory and Appl}. {\bf 45}, N. 4, 2000.


 \item{[S]} D. H. Sattinger,
  {\sl Monotone methods in Nonlinear Elliptic and Parabolic Boundary Value Problems.}.
 { Indiana  University Math. J. }. {\bf 21}, N. 11, 1972.

\item{[Sp]}  H. Spohn,  {\sl Large scale dynamics of interacting particles}.
(1991), Springer, Berlin.

\item{[Sp1]}   H. Spohn,  {\sl Long range correlations for stochastic lattice gases in a non-equilibrium steady state}.  {   J. Phys. A:Math. Gen.} {\bf16}  4275-4291,   1983  


\end